\documentclass[11pt]{amsart}

\usepackage{amsfonts,amssymb,amscd,amsmath,latexsym,amsbsy,amsthm,comment}

\usepackage{amssymb}
\usepackage{amsfonts}
\usepackage{latexsym}
\usepackage{epsfig}
\usepackage{xypic}

\theoremstyle{plain}
\newtheorem{theorem}{Theorem}
\newtheorem{thm}[theorem]{Theorem}

\newtheorem{lemma}[theorem]{Lemma}
\newtheorem{lem}[theorem]{Lemma}
\newtheorem{proposition}[theorem]{Proposition}
\newtheorem{prop}[theorem]{Proposition}

\newtheorem{corollary}[theorem]{Corollary}
\newtheorem{cor}[theorem]{Corollary}

\newtheorem{conjecture}[theorem]{Conjecture}

\theoremstyle{definition}
\newtheorem{example}[theorem]{Example}
\newtheorem{definition}[theorem]{Definition}

\theoremstyle{remark}
\newtheorem*{remark}{Remark}
\newtheorem*{remarks}{Remarks}

\numberwithin{theorem}{section}
\numberwithin{equation}{section}

\newcommand{\Tr}{\text{Tr}}

\newcommand{\rank}{\text{rank}}

\newcommand{\Hom}{\text{Hom}}

\newcommand{\g}{\mathfrak{g}}
\newcommand{\h}{\mathfrak{h}}

\newcommand{\ben}{\begin{enumerate}}
\newcommand{\een}{\end{enumerate}}

\newcommand{\QQ}{{\mathbb{Q}}}
\newcommand{\CC}{{\mathbb{C}}}
\newcommand{\RR}{{\mathbb{R}}}
\newcommand{\FF}{{\mathbb{F}}}
\newcommand{\ZZ}{{\mathbb{Z}}}
\newcommand{\RP}{{\mathbb{RP}}}
\newcommand{\CP}{{\mathbb{CP}}}

\def\Lambdasharp{\Lambda^{\!\scriptscriptstyle\#}}
\def\hGn{\widehat{\Gamma_n}}
\def\Piodd#1{\Pi^{\raisebox{.5ex}{\tiny\it odd}}_{#1}}

\newcommand{\solu}[1]{\begin{sol}{\bf (\ref{#1})}}

\begin{document}

\title{The cohomology ring of the real locus of the moduli space of stable curves of genus $0$
with marked points}

\author{Pavel Etingof}
\author{Andr\'e Henriques}
\author{Joel Kamnitzer}
\author{Eric Rains}

\maketitle
\tableofcontents

\section{Introduction}

Let $\overline{M_{0,n}}$ be the Deligne-Mumford compactification
of the moduli space of algebraic curves of genus 0 with $n$ labelled points. It is a smooth projective variety over $\Bbb Q$ (\cite{DM}) which parametrizes stable, possibly singular, curves of genus $0 $ with $ n $ labelled points.
The geometry of this variety is very well understood; in particular,
the cohomology ring of the manifold of complex points
$\overline{M_{0,n}}(\Bbb C)$ was computed by Keel \cite{Keel}, who also showed that
it coincides with the Chow ring of $\overline{M_{0,n}}$.

In this paper we will be interested in the topology of the
manifold $M_n:=\overline{M_{0,n}}(\Bbb R)$ of {\bf real} points of
this variety.  There are a number of prior results concerning the
topology of this manifold.  Kapranov and Devadoss (\cite{Kap},
\cite{Dev}) found a cell decomposition of $M_n$, and Devadoss used it to
determine its Euler characteristic.  Davis-Januszkiewicz-Scott
(\cite{DJS}) found a presentation of the fundamental group $\Gamma_n $ 
of $M_n$ and proved that $M_n$ is a $K(\pi, 1) $ space for this
group.  Nevertheless, the topology of $M_n$ is less well
understood than that of $\overline{M_{0,n}}(\Bbb C)$. In
particular,  the Betti numbers of
$M_n$ have been unknown until present.

\subsection{The cohomology of $M_n$}
In this paper, we completely determine the cohomology of $M_n$
with rational coefficients.  In Definition \ref{def:Lambda}, we give a presentation
of the cohomology algebra $H^*(M_n,\Bbb Q)$ by generators and
relations.  Later in Theorem \ref{th:basis}, we give a basis of this algebra and use this show that its
Poincar\'e series equals
$$
P_n(t)=\prod_{0\le k<(n-3)/2}(1+(n-3-2k)^2t).
$$
In particular, we show that the rank of $ H_1(M_n) $ is
$\binom{n-1}{3} $ which answers a question of Morava,
\cite[p.5]{Mo}.\footnote{Some partial results in this direction were
obtained in \cite{Ba1,Ba2}.}

The variety $ \overline{M_{0,n}}$, and hence the manifold $ M_n $, has an action of the symmetric group $ S_n $ which permutes the labelled points.   In a subsequent paper, E.R. \cite{R1} computes the character of the action of $S_n$
on $H^*(M_n,\Bbb Q)$.  

The manifolds $M_n$ are not orientable for $n\ge 5$, and in particular their cohomology
groups are not free and contain 2-torsion. We determine the 2-torsion in the cohomology,
and show there is no 4-torsion. The cohomology of $M_n$ does not have odd torsion (this has been recently shown by E.R. \cite{R2}),
so our results give a description of the cohomology $H^*(M_n,\Bbb Z)$.  
Note that since $M_n$ is a $K(\pi,1)$-space, this cohomology is also the cohomology of
the group $\Gamma_n$.

The description of the cohomology of $M_n $ has recently been generalized by E.R. \cite{R2} to a computation of the integral homology of the real points of any de Concini-Procesi model coming from any real subspace arrangement (the manifold $ M_n $ comes from the $A_{n-2} $ hyperplane arrangement).

\subsection{The operad structure}
The collection of spaces $M_n$ forms a topological operad, since
stable curves of genus $0$ can be attached to each other at marked
points, as described in section \ref{se:operad}.
Similarly, the homology $H_*(M_n,\Bbb Q)$ is an operad in the
symmetric monoidal category of $\Bbb Z$-graded $\Bbb
Q$-supervector spaces. This operad was first discussed by Morava
\cite{Mo}, who suggested that it might be related to symplectic
topology.  Understanding this operad was one of the primary
motivations for this work.

In Theorem \ref{th:operad}, we show that this operad generated by a
supercommutative associative product $ab$ of degree $0$ and
skew-supercommutative ternary ``2-bracket'' $[a,b,c]$, such that
the 2-bracket is a derivation in each variable and satisfies a
quadratic Jacobi identity in the space of 5-ary operations.  Motivated by the Hanlon-Wachs theory of Lie 2-algebras \cite{HW}, we call this the operad of 2-Gerstenhaber algebras. 

The structure of the
homology operad of $\overline{M_{0,n}}(\Bbb C)$ was determined by Konstevich and Manin (\cite{KM2}, see
also \cite{Ge}); in this case the operad (called the operad of
hypercommutative algebras) turns out to be infinitely generated.

\subsection{The analogy with braid groups}
We see that the space $\overline{M_{0,n}}(\Bbb R)$ has very
different topological properties from those of
$\overline{M_{0,n}}(\Bbb C)$. Indeed, $\overline{M_{0,n}}(\Bbb R)$
is $K(\pi,1)$, its Poincar\'e polynomial has a simple
factorization, its Betti numbers grow polynomially in $n$, and its
homology is a finitely generated operad.  In contrast,
$\overline{M_{0,n}}(\Bbb C)$ is simply connected, its Poincar\'e
polynomial does not have a simple factorization, its Betti numbers
grow exponentially, and its homology operad is infinitely
generated. 

On the other hand, the properties of the configuration
space $C_{n-1} = \mathbb{C}^{n-1}\smallsetminus \Delta $ of $n-1$-tuples of distinct complex numbers are
much more similar to those of $\overline{M_{0,n}}(\Bbb R)$.
It is a $K(\pi,1)$ space with Poincar\'e polynomial \linebreak
$(1+t)\cdots(1+(n-2)t)$, so its Betti numbers grow polynomially.  Also
its homology operad is well known to be the operad of
Gerstenhaber algebras, which has two binary generators. The
analogy between $M_n$ and $C_{n-1}$ and between their fundamental
groups (the pure cactus group $\Gamma_n$ and the pure braid group
$PB_{n-1}$), discussed already in \cite{Dev},\cite{Mo}, and
\cite{HK}, is very useful and has been a source of inspiration for
us while writing this paper.

\subsection{The Lie algebra $L_n$}
One application of the computation of the cohomology ring is that it allows us to understand various Lie algebras associated to the group $ \Gamma_n $.

The cohomology algebra $H^*(M_n,\Bbb Q) $ is a quadratic
algebra and thus we can consider its quadratic
dual algebra $U_n$, which is the universal enveloping algebra of a
quadratic Lie algebra $L_n$.  Since we have a presentation of $ \Lambda_n $, we get a presentation of $ L_n $ (see Proposition \ref{th:Liepresent}).

On the other hand, one can construct a Lie algebra $ \mathcal{L}_n $ directly from $ \Gamma_n $, by taking the associated graded of lower central series filtration and then quotienting by the 2-torsion.  In Theorem \ref{homlie}, we construct a surjective homomorphism of graded Lie
algebras  $\psi_n: L_n \to {\mathcal L}_n$. We expect that this
homomorphism is actually an isomorphism, similarly to the braid group case.  

We also expect that the algebra $U_n $ is Koszul.  On the other hand, somewhat disappointingly, we 
show that for $n\ge 6$ the Malcev Lie algebra of $\Gamma_n$ 
is not isomorphic to the degree completion of $L_n\otimes \Bbb Q$, 
and in particular the spaces $M_n$ for $n\ge 6$ are not formal. This fact reflects an essential 
difference between the pure cactus group and the pure braid group. 

\subsection{Relation to coboundary Lie 
quasibialgebras and quasiHopf algebras.}

The motivation for the conjecture that 
the map $\psi_n$ is an isomorphism comes from the theory of coboundary Lie quasibialgebras.
Let $\g$ be a Lie algebra over a field of characteristic zero, 
with a coboundary Lie quasibialgebra
structure $\varphi\in (\wedge^3\g)^\g$ (\cite{Dr1}).
Let $X_1, \cdots, X_{n-1}$ be representations of $\g$. 
From the explicit form 
of generators and relations for $L_n$, we show that $L_n$ acts 
on $X_1\otimes \cdots \otimes X_{n-1}$ (see section \ref{se:quasiquasi}). 

On the other hand, 
Drinfeld showed in \cite{Dr1} that any coboundary Lie
quasibialgebra can be quantized to a coboundary quasiHopf algebra.  The representation category of such a quasiHopf algebra is a coboundary category.  The group $ \Gamma_n $ acts on a tensor product 
$Y_1\otimes \dots \otimes Y_{n-1}$ in any coboundary monoidal
category (see \cite{HK} and section \ref{se:coboundary}).  From this, we get an action of $ \mathcal{L}_n $ on $ X_1 \otimes \dots \otimes X_{n-1} $.  

Theorem \ref{factor} shows that the action of $L_n$ on 
$X_1\otimes \dots \otimes X_{n-1}$ factors through the morphism
$\psi_n$ and this action of $ \mathcal{L}_n$. 
 
The above statements are direct analogs 
of the corresponding results for pure braid groups and quasitriangular Lie quasibialgebras, as developed by Drinfeld in \cite{Dr1,Dr2}.  Moreover, the entire action of the braid group on tensor products can be recovered as the monodromy of the Knizhnik-Zamolodchikov connection.  However, 
because of the non-formality of $M_n$, at the moment we are pessimistic about the existence of an analogous result in our case. 


\subsection{Organization}
 The paper is organized as follows. Section 2 contains
the statements of the main theorems describing
$H^*(M_n, \Bbb Q) $.  In section 3 we give additional results and conjectures, mostly concerning the Lie algebra $ L_n $.  Section 4, 5 and 6 are devoted to the proof of the main theorem.  First, in section 4, we prove that a certain algebra related to $H^*(M_n, \Bbb Q) $
has a basis indexed by
``basic triangle forests'' (combinatorial objects we introduce
for this purpose). In section 5, we recall Keel's description of 
$ H^*(\overline{M_{0,n}}(\Bbb C),\Bbb Z) $ and use it to
give an upper bound on the ranks of $ H^*(M_n, \Bbb Q)$.  Finally
in section 6, we prove our main results concerning the cohomology
ring of $M_n$. 

\subsection*{Acknowledgments} The authors are grateful to L.
Avramov, C. De Concini, J. Morava, J. Morgan, and B. Sturmfels, for useful
discussions and references. P.E. thanks the mathematics department
of ETH (Zurich) for hospitality. The work of P.E. was partially
supported by the NSF grant DMS-0504847 and the CRDF grant
RM1-2545-MO-03. E.R. was supported in part by NSF Grant No.
DMS-0401387. J.K. thanks the mathematics department of EPFL for
hospitality.  The work of J.K. was supported by NSERC and AIM.
Finally, we would like to mention that at many stages of this
work we made significant use of the Magma computer algebra system
for algebraic computations. 

\section{The cohomology ring and the homology operad}
\subsection{The algebra $\Lambda_n$}
We begin by introducing an algebra which will turn out to be
equal to the cohomology ring of $M_n$ over $ \mathbb{Q} $.

\begin{definition} \label{def:Lambda}
$\Lambda_n$ is the skew-commutative algebra
generated over $\Bbb Z$ by elements $\omega_{ijkl}$, $1\le i,j,k,l\le n$,
which are antisymmetric in $ijkl$, with defining relations
\begin{equation}\label{om1}
\omega_{ijkl}+\omega_{jklm}+\omega_{klmi}+\omega_{lmij}+\omega_{mijk}=0,
\end{equation}
\begin{equation}\label{om2}
\omega_{ijkl}\,\omega_{ijkm}=0,
\end{equation}
\begin{equation}\label{om3}
 \omega_{ijkl}\,\omega_{lmpi}
+\omega_{klmp}\,\omega_{pijk}
+\omega_{mpij}\,\omega_{jklm}
=
0
\end{equation}
for any distinct $i,j,k,l,m,p$,
\end{definition}

In particular, $\Lambda_n$ is quadratic. 

We will also consider the algebras $\Lambda_n\otimes R$ for commutative rings $R$.
They are defined over $R$ by the same generators and relations. 

\begin{remark} One can show (by a somewhat tedious calculation, which we
did using the program ``Magma'')
that 2 times (\ref{om3}) is in the ideal generated by
(\ref{om1}) and (\ref{om2}). So this relation becomes redundant if $1/2\in R$.
\end{remark}

The algebra $\Lambda_n$ has a natural action of $S_n$.

\begin{proposition}\label{is} One has
$\Lambda_n[1]=\wedge^3\h_n$, as $S_n$-modules, where $\h_n$ is the $n-1$-dimensional
submodule of the permutation representation, consisting of vectors with zero sum of coordinates
(in particular, $\Lambda_n[1]$ is free of rank $(n-1)(n-2)(n-3)/6$).
\end{proposition}

\begin{proof}
An isomorphism $\Lambda_n[1]\to \wedge^3\h_n$ is given by
$$
\omega_{ijkl}\to (e_i-e_l)\wedge(e_j-e_l)\wedge(e_k-e_l).
$$\end{proof}

We now switch to a different presentation of $\Lambda_n$.  In this presentation
only the $S_{n-1}$-symmetry, rather
than the full $S_n$-symmetry, is apparent.  However the presentation
only contains quadratic relations.

\begin{proposition}\label{nupres}
The algebra $\Lambda_n$ is isomorphic (in a natural way) to the
skew-commutative algebra generated by $\nu_{ijk}$, $1\le i,j,k\le
n-1$ (antisymmetric in $ijk$) with defining relations
\[
\nu_{ijk}\nu_{ijl} = 0,
\]
and
\[
 \nu_{ijk}\nu_{klm}
+\nu_{jkl}\nu_{lmi}
+\nu_{klm}\nu_{mij}
+\nu_{lmi}\nu_{ijk}
+\nu_{mij}\nu_{jkl}
=0.
\]
\end{proposition}

\begin{proof}
Let $\Lambda_n'$ be the algebra defined as
in the proposition. Define a homomorphism
$f: \Lambda_n'\to \Lambda_n$ by the formula
$\nu_{ijk}\mapsto \omega_{ijkn}$.
By directly manipulating the relations, it is easy to see 
that this homomorphism is well defined.  Using the 5-term linear relation (\ref{om1}), we can find an inverse for $ f $.
Thus $f$ is an isomorphism.
\end{proof}

\begin{theorem}\label{hilser}
For each $n$, $\Lambda_n$ is a free $\ZZ$-module with Poincar\'e polynomial
$$
P_n(t)=\prod_{0\le k<(n-3)/2}(1+(n-3-2k)^2t).
$$
\end{theorem}

The proof of this theorem is given in Section 6.

\subsection{The real moduli space}

\subsubsection{Stable curves}
Recall \cite{DM} that a stable curve of genus 0 with $n$ labeled points is a
finite union $C$ of projective lines $C_1,...,C_p$, together with labeled distinct points
$z_1,...,z_n\in C$ such that the following conditions are satisfied
\begin{enumerate}
\item each $z_i$ belongs to a unique $C_j$;

\item $C_i\cap C_j$ is either empty or consists of one point, and in the latter
case the intersection is transversal;

\item The graph of components (whose vertices are the lines $C_i $
and whose edges correspond to pairs of intersecting lines) is a
tree;

\item The total number of special points (i.e. marked points or
intersection points) that belong to a given component $C_i$ is at
least $3$.
\end{enumerate}

So a stable curve must have at least $3$ labeled points.

\begin{minipage}{12.5cm}
\vspace{.3cm}
\centerline{
\psfig{file=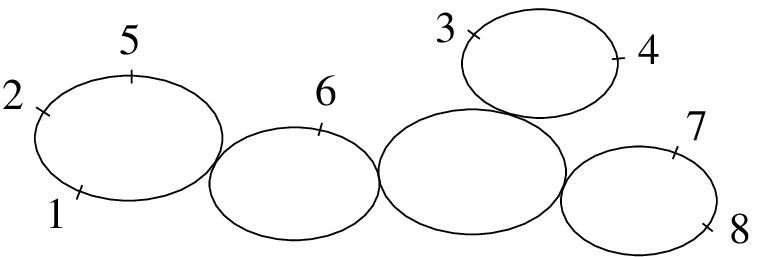,width=4cm}}

\centerline{{\it A stable curve with 8 marked points.}}
\vspace{.3cm}
\end{minipage}

An equivalence between two stable curves
$C=(C_1,...,C_p,z_1,...,z_n)$ and
$C'=(C_1',...,C_p',z_1',...,z_n')$ is an isomorphism of algebraic curves
$f: C\to C'$ which maps $z_i$ to $z_i'$ for each $i$. 
Thus $f$ reduces
to a collection of $p$ fractional linear maps $f_i: C_i\to
C_{\sigma(i)}'$, 
where $\sigma$ is a permutation.
It is easy to see that any equivalence of $C$ to itself is the identity.

Over the real numbers, 
the projective lines are circles, so a stable curve
is a ``cactus-like'' structure -- a tree of circles with labeled
points 
on them.

\subsubsection{The moduli space $M_n$}
For $n\ge 3$, let $M_n$  be the real locus
of the Deligne-Mumford compactification
of the moduli space of curves of genus zero
with $n$ marked points labeled $1,...,n$.
In other words, $M_n$ is the set of equivalence classes of stable
curves of genus 0 with $n$ labeled points defined over $\Bbb R$.
Clearly, $M_n$ carries a natural (non-free) action of $S_n$.

\begin{example}
\begin{enumerate}
\item $M_3$ is a point.

\item $M_4$ is a circle. More precisely, the cross-ratio map defines an
isomorphism $M_4\to \Bbb R\Bbb P^1$.

\item $M_5$ is a compact connected nonorientable surface
of Euler characteristic $-3$, i.e. the connected sum of 5 real
projective planes (see \cite{Dev}).
\end{enumerate}
\end{example}

The following theorem summarizes
some of the known results about $M_n$.

\begin{theorem}\label{known}
\begin{enumerate}
\item $M_n$ is a connected,
compact, smooth manifold
of dimension $n-3$.

\item  The Euler characteristic of $M_n$
is $0$ for even $n$ and
$$(-1)^{(n+1)/2}(n-2)!!(n-4)!!$$ if $n$ is odd.

\item $M_n$ is a $K(\pi,1)$-space.
\end{enumerate}
\end{theorem}

Part (i) is well-known and appears in \cite{Dev} and \cite{DJS}.
Part (ii) is due to Devadoss (see \cite[Theorem 3.2.3]{Dev}) and
Gaiffi (see \cite{Gai}) and comes from understanding the natural cell structure on $M_n $. Part (iii) is due to Davis-Januszkiewicz-Scott \cite{DJS}.  It is proven by showing that $ M_n $ is a Cat(0)-space.

\subsection{The cohomology of $M_n$}\label{sec223}
Let us now formulate the main result of this paper.
To do so, note that for any ordered $m$-element subset $S=\lbrace{s_1,...,s_m\rbrace}$
of $\lbrace{1,...,n\rbrace}$
we have a natural map $\phi_S: M_n\to M_m$, forgetting the points with
labels outside $S$. More precisely, given a stable curve $C$
with labeled points $z_1,...,z_n$, $\phi_S(C)$ is $C$ with labeled points
$z_{s_1},...,z_{s_m}$, in which the components that have fewer than 3 special points
have been collapsed in an obvious way.

Thus for any commutative ring $R$ we have a homomorphism of algebras
$\phi_S^*: H^*(M_m,R)\to H^*(M_n,R)$.

For $m=4$, $M_m$ is a circle, and we denote by $\omega_S$ the image of the
standard generator of $H^1(M_4,R)$ under $\phi_S^*$.

\begin{proposition}
Over any ring $R$ in which $2$ is invertible, the elements $\omega_S$ satisfy the relations
(\ref{om1}), (\ref{om2}).
\end{proposition}

\begin{proof} It is sufficient to consider the case $R=\Bbb Z[1/2]$.
The skew-symmetry of $\omega_S$ is obvious.

Next, we check the quadratic relations (\ref{om2}).   By
considering the maps $\phi_S$ for $|S|=5$, it suffices to check this
relation on $ M_5$. But $H^2(M_5,\Bbb Z)=\Bbb Z/2\Bbb Z$
because $M_5$ is non-orientable, so by the universal coefficient
theorem $H^2(M_5,R)=0$. 

The 5-term linear relation (\ref{om1})  may also be checked on
$M_5$. Since $H^1(M_5,\Bbb Z)$ is free over $\Bbb Z$, it is
sufficient to check the relation after tensoring with $\Bbb Q$. As an $S_5$-module, $H^1(M_5,\Bbb Q)$ is the
tensor product of the permutation and sign representations.  In
particular, the 5-cycle has no invariants in this representation,
and hence the 5-term relation holds.
\end{proof}

\begin{corollary} For any ring $R$ in which $2$ is invertible, we have a homomorphism of algebras
\begin{equation}\label{deffR}
f_n^R: \Lambda_n\otimes R\to H^*(M_n,R),
\end{equation}
which maps $\omega_S$ to $\omega_S$.
\end{corollary}

Our main result is the following theorem.

\begin{theorem}\label{cohom}
$f_n^\Bbb Q$ is an isomorphism.
\end{theorem}

It then follows from Theorem \ref{hilser} that the Poincar\'e polynomial of $M_n$
is $P_n(t)$.
In the process, we also prove the following result.

\begin{theorem}\label{cohom1}
$H^*(M_n,\Bbb Z)$ does not have $4$-torsion.
\end{theorem}

A description of the 2-torsion in $H^*(M_n,\Bbb Z)$ (which happens to be
quite big) is given in Section 5.

The absence of $4$-torsion in $H^*(M_n,\Bbb Z)$ allows us to
sharpen the above statements as follows.

For any abelian group $A$, we denote by $A\langle 2\rangle $ the 2-torsion in $A$.

\begin{proposition}
Over any ring $R$, the elements $\omega_S$ satisfy the relations
(\ref{om1}), (\ref{om2}), and (\ref{om3}) modulo 2-torsion.
Hence, we have a homomorphism of algebras
\[
f_n^R: \Lambda_n\otimes R\to H^*(M_n,R)/H^*(M_n,R)\langle 2\rangle ,
\]
which maps $\omega_S$ to $\omega_S$.
\end{proposition}

\begin{proof} It is sufficient to consider the case $R=\Bbb Z$.
Relations (\ref{om1}) and (\ref{om2}) modulo 2-torsion are proved in the
same way as in the case when $\frac{1}{2}\in R$.  It remains to prove
(\ref{om3}).  Although the relation (\ref{om3}) is not in the ideal
generated by the other relations (over $\Bbb Z$), if we multiply it by 2,
we obtain an element of this ideal, and thus the relation holds modulo
4-torsion.  Since by Theorem \ref{cohom1},
$H^*(M_6,\Bbb Z)$ has no 4-torsion\footnote{This can also be
verified directly, either from the chain complex of \cite{Dev} or using the
Bockstein map as discussed below.}, the result follows.
\end{proof}

Now Theorem \ref{cohom} can be strengthened as follows.
Let $\Bbb Z_{2}$ be the ring of 2-adic integers (we could 
equivalently consider the ring of 2-local integers).

\begin{theorem}\label{cohomm}
$f_n^{\Bbb Z_{2}}$ is an isomorphism.
\end{theorem}

In fact, one has the following stronger result. 

\begin{theorem}\label{odd} (\cite{R2}) $H^*(M_n,\Bbb Z)$ does not have odd
torsion. In particular, $f_n^{\Bbb Z}$ is an isomorphism.
\end{theorem}

\begin{remark} 
This theorem was conjectured in the first version of the present paper,
on the basis of a computation by A.H. and John Morgan, who
checked, using the blowup construction of $M_n$, that the theorem holds for
$n\le 8$.  This conjecture was recently proved by E.R. as 
a special case of the main theorem of \cite{R2}.
\end{remark}

\subsection{The operad structure on the homology of $M_n$.}
\label{se:operad}
The spaces $M_n$ form a topological operad which was first studied by Devadoss, who called it the mosaic operad
\cite{Dev}. To define this operad, it is convenient to agree that
each of the undefined moduli spaces $M_1$ and $M_2$ consists of
one point. We will define a topological operad with a space of
$n$-ary operations $M(n):=M_{n+1}$.  We think of a point of
$M_{n+1}$ as an $n$-ary operation where the inputs sit at points
$1,...,n$ and the output is $n+1$. 

The operad structure is defined by attaching curves at marked points. More explicitly, given $p,q\ge 0$ and
$1\le j\le p$, we have a ``substitution'' map $\gamma_{jpq}:
M_{p+1}\times M_{q+1}\to M_{p+q}$ given by attaching a curve $C_1$
with $p+1$ marked points to a curve $C_2$ with $q+1$ marked points
by identifying the point $j$ on the first curve with the point
$q+1$ on the second curve, and then adding $q-1$ to the labels $
j+1, \dots, p+1 $ on $C_1$ and adding $j-1 $ to the labels of the
points $1,...,q$ on $C_2$. The operad structure is obtained by
iterating such maps.

Recall that a cyclic operad (see \cite{GK}) is an operad $ P(\bullet) $ in which the action of $ S_n $ on $ P(n) $ extends to an action of $ S_{n+1} $, compatible with the operad structure.  In our case, $ S_{n+1} $ acts on $ M(n) = M_{n+1} $ in a natural way and thus $ M(\bullet) $ is a cyclic operad.

\begin{remark} For clarity, let us separately discuss special
cases of the substitution map when $p$ or $q\le 1$. If $q=1$, or $p=1$, $\gamma_{jpq}$ is
the identity map. If $q=0$, $\gamma_{jpq}: M_{p+1}\to M_p$ is the
map of erasing the $i$-th point.
\end{remark}

Since $M(\bullet)$ is a topological operad, the spaces
$O(n):=H_*(M_{n+1},\Bbb Q)=(\Lambda_{n+1}\otimes \Bbb Q)^*$ form an operad
in the category of $\Bbb Z$-graded supervector spaces. The
following result determines the structure of this operad.

\begin{theorem} \label{th:operad}
The operad $O(n)$ is the operad of unital
2-Gerstenhaber algebras. More specifically, it is generated
by $1\in O(0)$, $\mu\in O(2)$, and $\tau\in O(3)$, such that

\begin{enumerate}

\item $\mu$ is a commutative associative product of degree 0 with unit $1$;

\item $\tau$ is a skew-symmetric ternary operation
of degree $-1$, which is a derivation in
each variable with respect to the product $\mu$.

\item $\tau$ satisfies the Jacobi identity:
${\rm Alt}(\tau\circ (\tau\otimes {\rm Id}\otimes {\rm Id}))=0$,
where the alternator\footnote{As usual, the alternator is
understood in the supersense.} is over $S_5$.

\end{enumerate}
\end{theorem}

Algebras over this operad can be thought of as Lie 2-algebras with some additional structure.  More precisely, we have the following result.  Let $ \Sigma $ be the sign operad, which was introduced by Ginzburg-Kapranov \cite{GiK}.  
Let $ HW $ denote the Hanlon-Wachs operad of Lie 2-algebras (see \cite{HW}).  It it generated by a skew-symmetric ternary operadtion satisfying the Jacobi identity ${\rm Alt}(\tau\circ (\tau\otimes {\rm Id}\otimes {\rm Id}))=0$.

\begin{corollary} Consider the sub-operad $O'\subset O$,
with $O'(2k)=0$, \linebreak $O'(2k+1)=H_k(M_{2k+2},\Bbb Q)$. There is an isomorphism of operads $ O' = HW \otimes \Sigma $.
\end{corollary}

\begin{remark} The reason for tensoring with $ \Sigma $ is that Hanlon and Wachs consider an odd supersymmetric ternary
operation, while we need an odd
superalternating ternary operation.  As a result, the $S_n$-representation on $O'(n)$ is the tensor
product of the $S_n$ representation on the Hanlon-Wachs operad
$HW(n)$ with the sign representation.
However, the notions of $O'$-algebra and
$HW$-algebra are essentially equivalent:
one can go from one to the other by shifting by 1 the underlying $\ZZ$-graded super vector space.
Such a twisted version is
briefly mentioned in the beginning of \cite{HW}.
\end{remark}

\subsection{A few words on the proofs}
Theorems \ref{cohomm} and \ref{th:operad} are proved simultaneously in section 6.  The main idea is as follows, though for simplicity we give this outline over $ \QQ$.  First, in section \ref{se:basis}, we introduce a twisted version of the ring $ \Lambda_n $.  We find a basis for this ring and find its Hilbert series (Corollary \ref{cor:hilser_t}).  Then in section \ref{se:bound}, we use known results about the cohomology of the complex moduli space in order to give an upper bound on the Betti numbers of $ M_n $ (Corollary \ref{th:bound}).  In turns out that this upper bound exactly matches the dimensions of the graded pieces of $ \Lambda_n $.  Hence the map  $ f_n^\QQ : \Lambda_n \otimes \QQ \rightarrow H^*(M_n, \QQ) $ is between two graded rings with the dimensions of the left side bounding the dimensions of the right side.  So it suffices to prove that $ f_n $ is injective.  To do this, we use the operad structure on $ H_*(M_n, \QQ) $  to construct elements of $ H_*(M_n, \QQ) $ which pair ``upper-triangularily'' with the images under $ f_n $ of our basis for $ \Lambda_n $.  

\section{Results and conjectures concerning Lie algebras}
In this section we study the interplay between the cohomology ring of $ M_n $, its fundamental group $ \Gamma_n $, and theory of coboundary quasibialgebras.  Of main interest to us are different Lie algebras which one can associate to the manifold $ M_n$ and to its fundamental group $\Gamma_n $.  Much of what we do here is inspired by similar constructions of Kohno and Drinfeld in the configuration space/braid group setting.

\subsection{The quadratic dual to $\Lambda_n$ and the Lie algebra $ L_n $}
\begin{proposition}\label{th:Liepresent}
The quadratic dual $\Lambda_n^!$ of $\Lambda_n$ is the algebra
$U_n$ generated over $\ZZ$ by $\mu_{ijk}$, $1\le i,j,k\le n$,
which are antisymmetric in $ijk$, with defining relations
$$
\sum_{i}\mu_{ijk}=0, \
[\mu_{ijk},\mu_{pqr}]=0
$$
for distinct $i,j,k,p,q,r$ (with the obvious action of $S_n$).
It is also generated by $\mu_{ijk}$ with
$1\le i,j,k\le n-1$ with defining relations
\[
[\mu_{ijk},\mu_{pqi}+\mu_{pqj}+\mu_{pqk}]=0, \
[\mu_{ijk},\mu_{pqr}]=0.
\]
(with the obvious action of $S_{n-1}$ which extends to an action of $S_n$).
\end{proposition}

\begin{proof}
The equivalence of the first and second presentations of $U_n$
follows immediately by solving the linear relations for $\mu_{ijn}$.

To show that the algebra $U_n$ is dual to $\Lambda_n$,
it is convenient to use the $S_{n-1}$-invariant presentation of $\Lambda_n$,
which does not have linear relations. Recall that $ \Lambda_n[1] $ has basis $ \{ \mu_{ijk} \}$.  We let $ \{\nu_{ijk} \} $ denote the dual basis for $ \Lambda^!_n[1] = \Lambda_n[1]^* $.  By definition of quadratic dual, these $ \nu_{ijk} $ will generate $ \Lambda^!_n[1] $.  

To find the relations for $ \Lambda^!_n $, we must compute $ R^{\perp} \subset \Lambda^!_n[1] \otimes \Lambda^!_n[1] $, where $ R \subset \Lambda_n[1] \otimes \Lambda_n[1] $ are the relations for $ \Lambda_n $.  A convenient way to find $ R^{\perp} $, is to set $\Omega=\sum_{i<j<k}\mu_{ijk}\nu_{ijk}\in \Lambda_n^!\otimes \Lambda_n$.
Then the relations of $\Lambda_n^!$ are given by the formula
$[\Omega,\Omega]=0$. The above relations of $U_n$ are obtained
from this equation by a direct calculation.
\end{proof}

Let $L_n$ be the Lie algebra over $\Bbb Z$ generated by $\mu_{ijk}$ with
relations as above; we have $U_n=U(L_n)$.
Thus, $L_n\otimes \Bbb Q$ is the rational 
holonomy Lie algebra of $M_n$ in the sense of Chen,
see \cite{PS}. 

Recall that a $\Bbb Z_+$-graded algebra $A$ over
a field $k$ with $A[0]=k$ is called Koszul if ${\rm Ext}^i_A(k,k)$
(where $k$ is the augmentation module) sits in degree $i$ for all
$i\ge 1$. 

\begin{conjecture}\label{koz} The algebra $\Lambda_n\otimes \Bbb Q$
(or, equivalently, $U_n\otimes \Bbb Q$) is Koszul.
In particular, $U_n\otimes \Bbb Q$ has
Hilbert series
$$
P_n^!(t)=\frac{1}{P_n(-t)}=\prod_{0\le k<(n-3)/2}(1-(n-3-2k)^2t)^{-1}.
$$
\end{conjecture}
The second statement follows from the first by a general result about Koszul algebras.

\begin{remarks}
\begin{enumerate}
\item This conjecture is true for $n\le 6$, as in those degrees
$\Lambda_n$ has a quadratic Gr\"obner basis, so is Koszul.

\item The Hilbert series formula has been verified computationally in degree 3 for $n\le 9$.
\end{enumerate}
\end{remarks}

\begin{conjecture}\label{free}
$U_n$ and $L_n$ are free $\Bbb Z$-modules.
\end{conjecture}

\begin{remark} Let $A$ be a $\Bbb Z_+$-graded algebra over $\Bbb Z$ such that $A[i]$ are finitely generated free
$\Bbb Z$-modules for all $i$. One may define $A$ to be Koszul if
for each $j\ge 1$, $Ext^j_A(\Bbb Z,\Bbb Z)$ is a free $\Bbb
Z$-module living in degree $j$. This is equivalent to saying that
the algebras $A/pA$ are Koszul for all primes $p$. Thus we may
strengthen Conjectures \ref{koz} and \ref{free} by conjecturing
that $\Lambda_n$ and $U_n$ are Koszul over $\Bbb Z$.
\end{remark}

\subsection{Inductive nature of $ L_n$}
Computations using the ``Magma'' program suggest that the nonconstant
coefficients of the series
$P^!_n(t)/P^!_{n+1}(t)$ are all negative and that the virtual
character $1-P^!_n(g,t)/P^!_{n+1}(g,t)$ (where $P_n^!(g,t):=P_n(g,-t)^{-1}$) is actually a character in all
degrees. This suggests the following additional conjecture.

\begin{conjecture}\label{freeness}  
\begin{enumerate}
\item The kernel of the natural morphism 
$U_{n+1}\to U_n$ (sending $\mu_{ijn}$ to zero
and $\mu_{ijk}$ to themselves for $i,j,k<n$) is a
free $U_{n+1}$-module.  

\item The kernel of the natural morphism
$L_{n+1}\to L_n$ is a free Lie algebra (with infinitely many generators).
\end{enumerate}
\end{conjecture}

We note that the two statements of the conjecture become equivalent 
after extension of scalars from $\Bbb Z$ to $\Bbb Q$. 
This follows from the following lemma. 

\begin{lemma} Let $L$ be a $\Bbb N$-graded Lie algebra
acting on another $\Bbb N$-graded Lie algebra $F$, 
both having finite dimensional graded pieces, 
and let $L\ltimes F$ be their semidirect product. 
Then the following conditions are equivalent. 

(i) $F$ is a free Lie algebra.

(ii) The kernel $K$ of the natural map $U(L\ltimes F)\to U(L)$ 
is a free $U(L\ltimes F)$-module. 
\end{lemma}

Indeed, the equivalence follows by applying the Lemma 
for $L=L_n$, $F=\ker(L_{n+1}\to L_n)$ and $L\ltimes F=L_{n+1}$. 

\begin{proof} We have a natural isomorphism 
$U(L\ltimes F)=U(L)\ltimes U(F)$,      
under which the kernel $K$ is identified with 
$U(L)\otimes U(F)F$.

If $F$ is a free Lie algebra,
 then it is freely generated by a graded subspace $G$. Hence 
$U(F)F=U(F)G$ is a free $U(F)$-module generated by $G$.
It follows that $K$ is a free $U(L\ltimes F)$-module generated by
$G$. Thus (i) implies (ii). 

Conversely, assume that $K$ is a free module over $U(L\ltimes F)$. 
It is clear that $F+(L\ltimes F)K=K$. 
Let $G\subset F$ be a graded complement to $(L\ltimes F)K$ in $K$. 
Then $K$ is freely generated by $G$, i.e., $K=U(L\ltimes F)\otimes
G=U(L)\otimes U(F)\otimes G$. Therefore $U(F)F$ is generated by
$G$ over $U(F)$. This implies that $F$ is generated by 
$G$ as a Lie algebra, and hence (by the Hilbert series
consideration) that $F$ is freely generated by $G$.
Thus (ii) implies (i).
\end{proof}

\begin{remark} Note that since $P^!_n(t)/P^!_{n+2}(t)=1-(n-1)^2t$,
it is tempting to make a much more simple-looking conjecture,
namely that the kernel of the homomorphism $L_{n+2}\to L_n$ is a
free Lie algebra generated by $(n-1)^2$ generators in degree $1$.
This, unfortunately, is very far from being true, since the graded
$S_n$-character of the kernel of this homomorphism is not the same
as that of a free Lie algebra (and in particular
$1-P^!_n(g,t)/P^!_{n+2}(g,t)$ is a virtual character, which is not
a character, and is not concentrated in degree $1$).
\end{remark}

\subsection{The rational $K(\pi,1)$-property}

A connected topological space $X$ is said to be rational $K(\pi,1)$
if its $\Bbb Q$-completion is a $K(\pi,1)$-space (\cite{BK}).
Note that a $K(\pi,1)$-space in the usual sense may not be a rational $K(\pi,1)$ space
(e.g., the complement of the complex hyperplane arrangement of type $D_n$).

It is proved in \cite{PY}, Proposition 5.2,
that for a connected topological space $X$ with finite Betti numbers,
if $H^*(X,\Bbb Q)$ is a Koszul algebra then $X$ is a rational $K(\pi,1)$ space.
Thus Conjecture \ref{koz} implies

\begin{conjecture}\label{rkp1}
The space $M_n$ is rational $K(\pi,1)$ in the sense of \cite{BK}.
\end{conjecture}

It is also shown in \cite{PY} that the Koszul property
of $H^*(X,\Bbb Q)$ and the rational $K(\pi,1)$ property of $X$
are equivalent if $X$ is a formal space. However, as we show below, 
the spaces $M_n$ are not formal for $n\ge 6$. 

\subsection{The fundamental group $\Gamma_n$ of $M_n$ and coboundary categories} \label{se:coboundary}

Let $\Gamma_n$ be the fundamental group of $M_n$.
To understand this group, we consider
another group $J_n$ which is the orbifold fundamental group
of the orbifold $M_{n+1}/S_n$ (the group $S_n$ leaves the point $n+1$ fixed).
There is a short exact sequence
$$
1\to \Gamma_{n+1}\to J_n\to S_n\to 1.
$$
Furthermore, it is explained in \cite{Dev,DJS,HK} that the group
$J_n$ has the following presentation: it is generated by elements
$s_{p,q}$, $1\le p<q\le n$, with defining relations
\begin{enumerate}
\item $s_{p,q}^2=1$;
\item $s_{p,q}s_{m,r}=s_{m,r}s_{p,q}$ if $[p,q]\cap [m,r]=\emptyset$;
\item $s_{p,q}s_{m,r}=s_{p+q-r,p+q-m}s_{p,q}$ if $[m,r]\subset [p,q]$.
\end{enumerate}

The above map $J_n\to S_n$ is
defined by sending $s_{p,q}$ to the involution that reverses the interval $[p,q]$ and keeps
the indices outside of this interval fixed.  The group $ J_n $ is called the ``cactus group'' and it is analogous to the braid group.

One significance of this group $ J_n $ comes from the theory of coboundary monoidal categories.  Recall (\cite{Dr1}, see also \cite{HK}) that a coboundary
monoidal category is a monoidal category $\mathcal C$ together
with a commutor morphism $c_{X,Y}: X\otimes Y\to Y\otimes X$,
functorial in $X,Y$, such that $c_{X,Y}\,c_{Y,X}=1$, and
$$
c_{Y\otimes X,Z}\,c_{X,Y}=c_{X,Z\otimes Y}\,c_{Y,Z}.
$$
(for simplicity we drop the associativity isomorphisms, assuming
that the category is strict, and write $c_{X,Y},c_{Y,Z}$ instead
of $c_{X,Y}\otimes 1_Z$, $1_X\otimes c_{Y,Z}$).

Let $ \mathcal C $ be such a category, and let $ X_1, \dots, X_n $
be $n $ objects in $ \mathcal C $.  Then, as shown in \cite{HK},
every element $g$ of the group $J_n$ defines a morphism
$X_1 \otimes \cdots \otimes X_n\to X_{g(1)}\otimes
\cdots \otimes X_{g(n)}$ (here $g(i)$ is the action of the image
of $g$ in $S_n$ on the index $i$). Namely, $s_{p,q}$ acts by
$$
c_{X_p,X_{p+1}}c_{X_p\otimes X_{p+1}, X_{p+2}} \cdots
\,c_{X_p\otimes \cdots \otimes X_{q-1},X_q}.
$$
The action of the cactus group on tensor products in coboundary monoidal categories is analogous to the action of the braid group on tensor products in braided monoidal categories.

For $1\le p\le q<r\le n$, let
$\sigma_{p,q,r}=s_{p,r}s_{p,q}s_{q+1,r}$ (we agree that
$s_{pp}=1$). Clearly, such elements generate $J_n$. 
The element $\sigma_{p,q,r}$ acts on $X_1\otimes \cdots \otimes X_n$
by $c_{X_p\otimes \cdots \otimes X_q,X_{q+1}\otimes \cdots \otimes
X_r}$, and its inverse is given by
$\sigma_{p,q,r}^{-1}=\sigma_{p,p+r-1-q,r}$.

Now for $1\le p\le q<r<m\le n$,
define
$$
b_{p,q,r,m} =
\sigma_{p,q,r}^{-1}
\sigma_{p+r-q,r,m}^{-1}
\sigma_{p,q,m}
$$
It is easy to see that $b_{p,q,r,m}\in \Gamma_n$. It acts on the
tensor product $X_1\otimes\cdots\otimes X_n$ by the morphism
$$
c_{Z,Y}\,c_{T,Y}\,c_{Y,Z\otimes T},
$$
where $Y=X_p\otimes\cdots\otimes X_q$,
$Z=X_{q+1}\otimes\cdots\otimes X_r$,
$T=X_{r+1}\otimes\cdots\otimes X_m$.
More geometrically, $b_{p,q,r,m}$ is represented by the following loop in $M_n$:

\begin{equation}\label{pictureloop}
\begin{matrix}\psfig{file=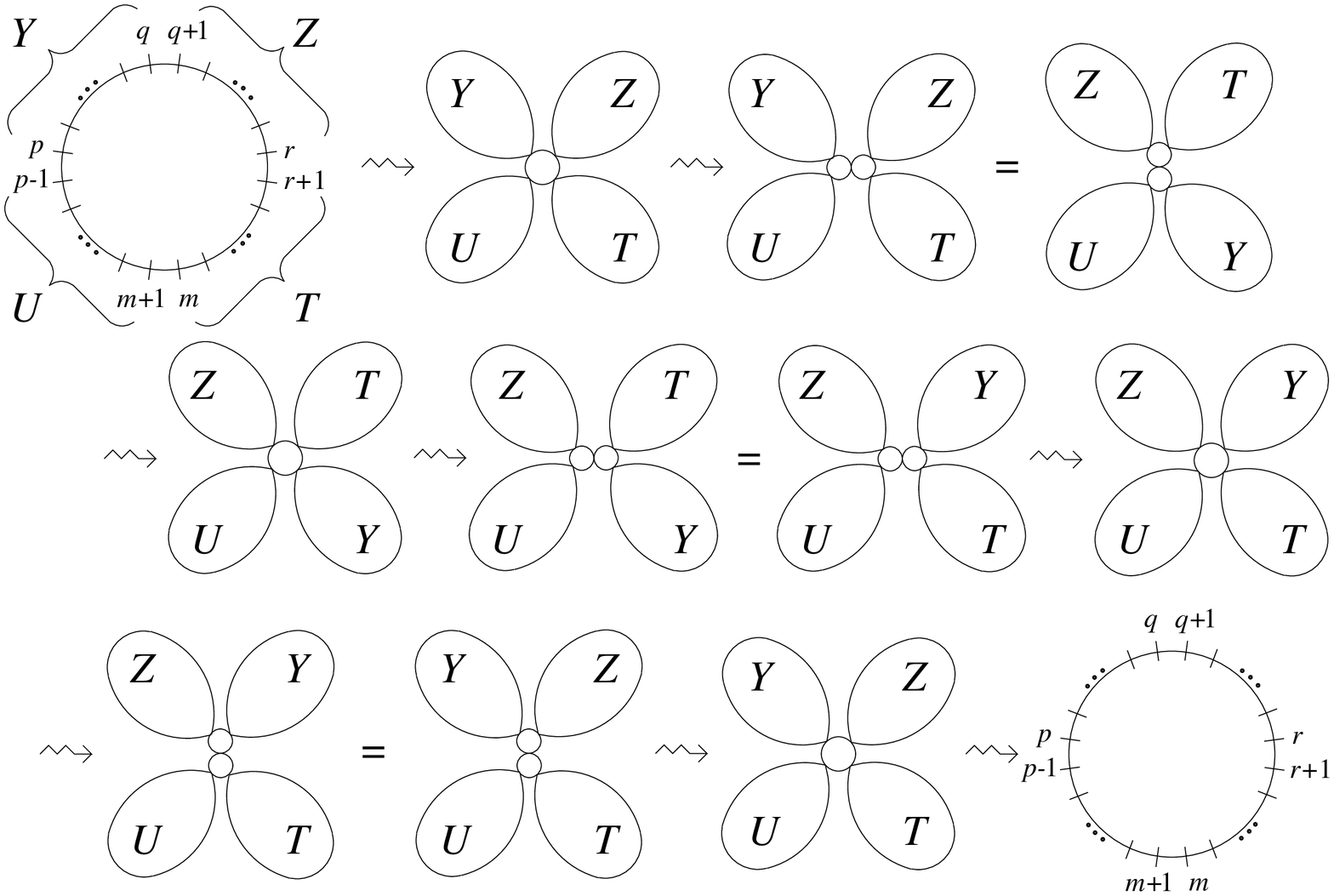,width=10cm}\end{matrix}\phantom{****}
\end{equation}

The significance of these elements is the following result.

\begin{prop} \label{th:genGamma}
The conjugates of $ b_{p,q,r,m} $ generate $ \Gamma_n $.
\end{prop}

\begin{proof}
Consider the quotient $G$ of $J_{n-1}$ by the relations
$b_{p,q,r,m}=1$.  We claim that $ G $ is isomorphic to $ S_{n-1} $. To show this, note that
$\sigma_{p,q,m}$ for $m>p+1$ can be expressed in $G$ via elements
$\sigma_{p',q',m'}$ with $m'-p'<m-p$. Thus the group $G$ is generated by $\sigma_{p,p,p+1}=s_{p,p+1}$,
which are easily shown to satisfy the braid relations, as desired.  Hence $ \Gamma_n $ is generated by the conjugates of $ b_{p,q,r,m} $.
\end{proof}

\subsection{The first homology of $M_n$}

\begin{theorem}\label{firsthom}
One has $H_1(M_n,\Bbb Z)=\wedge^3\Bbb Z^{n-1}\oplus E_n$, where $E_n$ is a vector space over $\Bbb F_2$.
\end{theorem}

\begin{proof} For every triple $i,j,k$ of distinct
indices between $1$ and $n-1$, we have
a submanifold of $D_{ijk}$ of $M_n$
(of codimension 1) which is the closure
of the set of curves with two components, one containing the points $i,j,k$ and
the other containing the remaining labeled points.
This submanifold is naturally isomorphic to $M_4\times M_{n-2}=S^1\times M_{n-2}$.
The circle $S^1$ has a natural orientation, which depends on the cyclic ordering of
$i,j,k$. Let $\mu_{ijk}$ denote the image in $H_1(M_n,\Bbb Z)$
of the fundamental class of this circle.
It is easy to see that $\mu_{ijk}$ are skew-symmetric in the indices
and are permuted in an obvious way by $S_{n-1}$.
Let $\nu_{ijk}\in H^1(M_n,\mathbb Q)$ be the images of $\nu_{ijk}$ under the map (\ref{deffR}).
Recall the map $\phi_S$ from section \ref{sec223}. The composite 
$$
S^1\hookrightarrow S^1\times M_{n-2}=D_{ijk}\hookrightarrow M_n
\stackrel{\phi_{i'\!j'\!k'\!n}}{-\!\!\!-\!\!\!\twoheadrightarrow}M_4=S^1
$$ 
is a diffeomorphism if $\{i,j,k\}=\{i',j',k'\}$ and is constant otherwise.
It follows that $\langle\mu_S,\nu_{S'}\rangle=\pm\delta_{S,S'}$ and that the $\mu_{ijk}$ are linearly independent over $\mathbb Q$.
It remains to show that the
$\mu_{ijk}$ for $1\le i,j,k<n$ span $H_1(M_n,\Bbb Z)/H_1(M_n,\Bbb Z)\langle 2\rangle $ (note that they do not span
$H_1(M_n,\Bbb Z)$ which has a very big 2-torsion).
To do so we will use the structure of the fundamental group $\Gamma_n$.

By Proposition \ref{th:genGamma}, the images
of the $b_{p,q,r,m}$ in $\Gamma_n/[\Gamma_n,\Gamma_n]$ generate $ \Gamma_n/[\Gamma_n,\Gamma_n]$.  So it suffices to show that each $ [b_{p,q,r,m}] $ can be written as a linear combination of the $ \mu_{ijk}$.  We claim that
\begin{equation}\label{sumf}
[b_{p,q,r,m}]=\sum_{p\le i\le q<j\le r<k\le m}\mu_{ijk}
\end{equation}
modulo 2-torsion.

If $n=4$, then $\Gamma_n=\mathbb Z$ and the
equality is exact (with only one term in the sum).
If $n=5$, then $H_1(M_n,\Bbb Z)=\Bbb Z^4\oplus \Bbb F_2$ (since $M_5$ is the
connected sum of 5 real projective planes)
and it's therefore enough to check (\ref{sumf}) over $\mathbb Q$.
To do so, we compute $\langle b_{p,q,r,m},\nu_{ijk}\rangle$.
Namely, we compose (\ref{pictureloop}) with $\phi_{ijk5}:M_5\to M_4$ and see that it produces
$1\in H_1(M_4)$ if $p\le i\le q<j\le r<k\le m$ and zero otherwise.

Consider now the case $n>5$.
It is easy to prove the following equalities:
\[
\begin{split}
b_{p,q,r,m}^{-1}&=\sigma_{p,q,r}^{-1}b_{p,p+r-q-1,r,m}\sigma_{p,q,r}\\
&=\sigma_{q+1,r,m}^{-1}b_{p,q,m+q-r-1,m}\sigma_{q+1,r,m}.
\end{split}
\]
Also, using the operad structure of $M_n$ and the case
$n=5$, we get that for any $p\le \ell< q$,
$$
[b_{p,q,r,m}]=[b_{\ell+1,q,r,m}]+\sigma_{p,\ell,q}^{-1}([b_{p+q-\ell,q,r,m}])
$$
in $\Gamma_n/[\Gamma_n,\Gamma_n]$ modulo $2$-torsion.
These two equalities imply equation (\ref{sumf}) by induction on
$m-p$. The theorem is proved.
\end{proof}

In particular, this gives a proof of Theorem 
\ref{odd} for $H_1(M_n,\Bbb Z)$.

\subsection{The lower central series of $\Gamma_n$}

Let $\Gamma_n^p$ be the $p$-th term of the lower central series of $\Gamma_n$,
and $L_n'$ be the associated graded $\Bbb Z$-Lie algebra of this series, i.e.
$L_n'=\oplus_{p=1}^\infty \Gamma_n^p/\Gamma_n^{p+1}$.
This is a Lie algebra graded by the positive integers,
and it is generated in degree $1$ by definition. 
Also we have a natural action
of $S_{n-1}$ on $L_n'$, coming from the group $J_{n-1}$.

We have $L_n'[1]=H_1(M_n,\Bbb Z)$, so $L_n'$ has 2-torsion. Let us
therefore consider the quotient Lie algebra ${\mathcal
L}_n=L_n'/L_n'\langle 2\rangle$. It is generated by its degree $1$
part, which by Theorem \ref{firsthom} is a free $\Bbb Z$-module
with basis $\mu_{ijk}$, $i<j<k$.

\begin{theorem}\label{homlie} There is a surjective $S_{n-1}$-equivariant homomorphism
of graded Lie algebras $\psi_n: L_n\to {\mathcal L}_n$, which maps
$\mu_{ijk}$ to $\mu_{ijk}$.
\end{theorem}

\begin{proof}
The only thing we need to prove is that $\mu_{ijk}\in {\mathcal L}_n[1]$ satisfy
the quadratic relations of $L_n$. Because of the $S_{n-1}$-symmetry, it is sufficient to show that
\begin{equation}\label{12345}
[\mu_{123},\mu_{145}+\mu_{245}+\mu_{345}]=0,\quad [\mu_{123},\mu_{456}]=0.
\end{equation}
To do so, note that, as was explained in the proof of Theorem \ref{firsthom},
$\mu_{123},\mu_{456}$ are the classes in
${\mathcal L}_n[1]$ of the elements $b_{1,1,2,3}$, $b_{4,4,5,6}$, while
$\mu_{145}+\mu_{245}+\mu_{345}$ is the class of the element
$b_{1,3,4,5}$. Now the identities (\ref{12345}) follow
from the relations $b_{1,1,2,3}b_{1,3,4,5}=b_{1,3,4,5}b_{1,1,2,3}$ and $b_{1,1,2,3}b_{4,4,5,6}=
b_{4,4,5,6}b_{1,1,2,3}$ in $\Gamma_n$, which easily follow from
the defining relations of $J_{n-1}$.
\end{proof}

\begin{conjecture}\label{isomo} \begin{enumerate}
\item $\psi_n$ is an isomorphism.
In particular, (assuming Conjecture \ref{free}),
${\mathcal L}_n$ is a free $\Bbb Z$-module, and thus
the only torsion in the lower central series of $\Gamma_n$ is 2-torsion.

\item $\cap_{k\ge 1}\Gamma_n^k=\lbrace{1\rbrace}$. In other words, the group $\Gamma_n$ is residually nilpotent.
\end{enumerate}
\end{conjecture}

\subsection{The prounipotent completion}

Let $\widehat{\Gamma_n}$ be the prounipotent 
(=Malcev) completion of $\Gamma_n$ (over $\Bbb Q$).
This is a prounipotent proalgebraic group. 
We can also define the proalgebraic group 
$\widehat{J_{n-1}}=J_{n-1}\times_{\Gamma_n}\widehat{\Gamma_n}$. 

Let ${\rm Lie}\widehat{\Gamma_n}$ be the Lie algebra
of the group $\widehat{\Gamma_n}$, which is called the Malcev Lie algebra of $\Gamma_n$ (see \cite{PS} for a more general discussion about such Lie algebras).  Let ${\rm grLie}\widehat{\Gamma_n}$ be the associated graded 
of this Lie algebra
with respect to the lower central series filtration.
It is easy to see that 
${\rm grLie}\widehat{\Gamma_n}={\mathcal L_n}\otimes \Bbb Q$.
Thus we have a surjective homomorphism $\psi_n^{\Bbb Q}: L_n\otimes \Bbb Q\to
{\rm grLie}\widehat{\Gamma_n}$, and we can make the $\Bbb Q$-version of
Conjecture \ref{isomo}:

\begin{conjecture}\label{isomo1} \begin{enumerate}
\item $\psi_n^{\Bbb Q}$ is an isomorphism.
\item $\Gamma_n$ injects into its prounipotent completion.
\end{enumerate}
\end{conjecture}

\begin{remark} Note that since $\Gamma_n$ is the fundamental group of 
an aspherical manifold, it is torsion free. This is a necessary condition for 
part (ii) of the conjecture.  
\end{remark}

\subsection{Coboundary Lie quasibialgebras, quasiHopf algebras
and representations of $L_n$ and
$\Gamma_n$.} \label{se:quasiquasi}

Recall \cite{Dr1} that a coboundary Lie quasibialgebra is a Lie
algebra $\g$ 
together with an element\footnote{We consider Lie quasibialgebras up to
twisting.} $\varphi\in
(\wedge^3\g)^\g$.  Given such a $(\g,\varphi)$
(over a field of characteristic zero), one can define a family of
homomorphisms $\beta_{n,\g,\varphi}: L_n\to U(\g)^{\otimes
{n-1}}$, defined by the formula
$\beta_{n,\g,\varphi}(\mu_{ijk})=\varphi_{ijk}$. (Here $ \varphi_{ijk} $ denotes the image of $ \varphi $ under the embedding $ \g^{\otimes 3} \rightarrow U(\g)^{\otimes n-1}$ which puts 1s in all factors other than $ i,j,k$.)  The invariance of $\varphi$ implies that 
the quadratic relations of $L_n$ are preserved under
this assignment.

\begin{theorem}\label{factor} The representations $\beta_{n,\g,\varphi}$ factor
through ${\mathcal L}_n$.
\end{theorem}

\begin{proof} The proof is based on the fact, due to Drinfeld \cite{Dr1}, that the representations
$\beta_{n,\g,\varphi}$ can be quantized, by quantizing the Lie
quasibialgebra $(\g,\varphi)$. Namely, Drinfeld showed that there
exists an associator
$\Phi=\Phi(\hbar^2)=1+\hbar^2\varphi/3+O(\hbar^4)$ in
$U(\g)^{\otimes 3}[[\hbar]]$ (given by some universal
formula in terms of $\hbar^2\varphi$), such that $(U(\g)^{\otimes
3}[[\hbar]],\Phi)$ is a coboundary quasiHopf algebra. (See
\cite{Dr1}, Proposition 3.10). 

Let $\mathcal C$ be the associated
category (the objects are representations of $\g$ and the
morphisms are power series in $\hbar$ whose coefficients are
morphisms of representations). This is a coboundary category and hence as explained in section \ref{se:coboundary}, there is an action of $ \Gamma_n $ on tensor products of objects.  In fact, from Drinfeld's construction, this action comes from a map $B_{n,\g,\varphi}: \Gamma_n\to 1+\hbar^2
U(\g)^{\otimes n-1}[[\hbar]]$. This homomorphism factors through
the prounipotent completion $\hGn$ of $\Gamma_n$, and
thus defines a Lie algebra homomorphism $B_{n,\g,\varphi}: {\rm
Lie}\hGn\to \hbar^2U(\g)^{\otimes n-1}[[\hbar]]$.

Under this homomorphism, the image of the $p$-th
term ${\rm Lie}\hGn^p$ of the lower central series
of ${\rm Lie}\hGn$ is contained in
$\hbar^{2p}U(\g)^{\otimes n-1}[[\hbar]]$. Therefore, we have a
natural homomorphism $B_{n,\g,\varphi}^0: {\rm grLie}\hGn
\to U(\g)^{\otimes n-1}$, defined by the formula
$B_{n,\g,\varphi}^0(z)=\hbar^{-2p}B_{n,\g,\varphi}(\hat z) {\rm
mod}\hbar$, where $z$ is of degree $p$ and $\hat z$ is a lift of
$z$. 

Now we claim that 
$\beta_{n,\g,\varphi}=B_{n,\g,\varphi}^0\circ \psi_n$. Both sides 
are invariant under the $S_{n-1}$-action, and thus it suffices to 
check the equality on the element $\mu_{123}$. This is straightforward to do by computing 
the $\hbar^2$-part of $B_{n,\g,\varphi}(b_{1,1,2,3})$ using that 
$\Phi=1+\hbar^2\varphi/3+...$. 
\end{proof}

\begin{remark} This theorem gives supporting evidence for the above conjecture
that $\psi_n$ is injective.
\end{remark}

\subsection{The non-formality of $M_n$.}

Let $\widehat{L_n\otimes \Bbb Q}$ be the degree completion of $L_n\otimes \Bbb Q$.

\begin{proposition} For $n\ge 6$, there does not exist 
an isomorphism \linebreak $\xi: {\rm Lie}\hGn
\to \widehat{L_n\otimes \Bbb Q}$  
whose associated graded is the identity in degree $1$. 
\end{proposition}

\begin{proof}
It suffices to prove the statement for $n=6$, 
since we have homomorphisms $\Gamma_{n-1}\to \Gamma_n$ and 
$\Gamma_n\to \Gamma_{n-1}$ whose composition is the identity
map $\Gamma_{n-1}\to \Gamma_{n-1}$, and similar homomorphisms for $L_n$.  

For $n=6$, the statement was checked computationally using 
the ``Magma'' program. 

More precisely, consider the exact sequence
$$
1\to \widehat{\Gamma_6}\to \widehat{J_5}\to S_5\to 1.
$$
Since $\widehat{\Gamma_6}$ is prounipotent over $\Bbb Q$ and $S_5$ is finite, 
this exact sequence is split. Thus 
we get 
$\widehat{J_5}=S_5\ltimes \widehat{\Gamma_6}$, for some action of $S_5$ 
on $\widehat{\Gamma_6}$.

Now assume that an isomorphism $\xi$ exists. Then 
all such isomorphisms form a torsor over a prounipotent
group of unipotent automorphisms of the target, 
so by averaging we can choose $\xi$ to be $S_5$-invariant.
Then $\xi$ can be lifted to a homomorphism 
$\xi':\widehat{J}_5\to S_5\ltimes  \exp(\widehat{L_6\otimes \Bbb Q})$.
This can be restricted to a homomorphism 
$\xi': J_5\to S_5\ltimes  \exp(\widehat{L_6\otimes \Bbb Q})$,
which maps $b_{1,1,2,3}$ into $\exp(\mu_{123}+...)$,
where $...$ are higher order terms. 
Since $J_5$ is given by simple generators and relations, 
one is able to look 
for such homomorphisms using ``Magma''
(modulo $m$-th commutators for some $m$, to make the computations finite). 
The computation showed that already modulo the fourth commutators 
there is no such homomorphism 
under which the element 
$b_{1,1,2,3}$ goes to $\exp(\mu_{123}+\text{higher terms})$.
This implies the required statement. 
\end{proof}

In particular, the proposition implies that for $n\ge 6$ 
the Malcev Lie algebra ${\rm Lie}\hGn$
of $\Gamma_n$ is not isomorphic to 
the degree completion of the 
rational holonomy Lie algebra $L_n\otimes {\Bbb Q}$ of $M_n$, which 
means that the group $\Gamma_n$ is not 1-formal in the sense of \cite{PS}.
Therefore, by the results of \cite{M,Su}, 
$M_n$ is not a formal space in the sense of Sullivan
(see e.g. \cite{PS} for more explanations). 

\begin{remark} For $n\le 5$, the spaces $M_n$ are formal
and the Malcev Lie algebra of $\Gamma_n$ is isomorphic to its 
associated graded. 
\end{remark}

\subsection{The Drinfeld-Kohno Lie algebra}

The Drinfeld-Kohno Lie algebra $\mathbf L_n$, defined for $n>2$, is
generated over $\Bbb Z$ by generators
$t_{ij}=t_{ji}$ for distinct indices $1\le i,j\le n-1$, and relations
$$
[t_{ij},t_{ik}+t_{jk}]=0,\ [t_{ij},t_{kl}]=0
$$
for distinct $i,j,k,l$ (see \cite{Ko,Dr1,Dr2}).
${\bold L}_n$ is a free $\Bbb Z$-module (since it is an iterated semidirect product of free
Lie algebras).

Kohno proved in \cite{Ko} (using the Knizhnik-Zamolodchikov
equations) that $\mathbf L_{n+1}\otimes \Bbb Q$ is the Lie algebra of the
prounipotent completion of the pure braid group $PB_n$. Let $\bold
U_n$ be the universal enveloping algebra of $\bold L_n$. The
algebra $\mathbf U_n\otimes \Bbb Q$ is Koszul, and is the
quadratic dual of the Orlik-Solomon algebra $OS_{n-1}$, which is
the cohomology algebra of the configuration space
$C_{n-1}=\lbrace{(z_1,...,z_{n-1})\in \CC^{n-1}: z_i\ne z_j\rbrace}$.
The Hilbert series of $\mathbf U_n$ is
$$
\mathbf P_n(t)=\prod_{m=1}^{n-2}(1-mt)^{-1}.
$$
(see e.g. \cite{Yuz2}.)

If $\g$ is a Lie algebra and $\Omega\in (S^2\g)^\g$ is an
invariant element (i.e., $(\g,\Omega)$ is a quasitriangular Lie quasibialgebra \cite{Dr1}),
then we have a homomorphism $\gamma_{n,\g,\Omega}: \mathbf L_n\to
U(\g)^{\otimes n-1}$ given by
$\gamma_{n,\g,\Omega}(t_{ij})=\Omega_{ij}$ (this is analogous to $ \beta_{n, \g, \varphi} $ above).  In the braid group setting, as $ \mathbf{L}_n $ is isomorphic to the Lie algebra of the prounipotent completion of $ PB_n$, the analogous statement to Theorem \ref{factor} is that the two different representations of $ \mathbf{L}_n $ agree (one comes from $\gamma$, the other from quantizing $\Omega$).

\subsection{Relation between $L_n $ and $\mathbf{L}_n$}
We have a homomorphism of Lie algebras
$\xi_n: L_n\to \mathbf L_n$ given by the formula
$\xi_n(\mu_{ijk})=[t_{ij},t_{jk}]$.
Note that under this homomorphism, $L_n[1]$ is identified with $\mathbf L_n[2]$.
It is clear that $\gamma_{n,\g,\omega}\circ \xi_n=\beta_{n,\g,\varphi}$,
where $\varphi=[\Omega^{12},\Omega^{23}]$.

\begin{theorem}\label{facto} The map
$\xi_n^{\Bbb Q}: L_n\otimes \Bbb Q\to {\mathbf L}_n\otimes \Bbb Q$
factors through ${\mathcal L}_n\otimes \Bbb Q$.
\end{theorem}

\begin{proof} Using Drinfeld's ``unitarization'' trick, 
one can define a homomorphism of
prounipotent completions $\Xi_n: \hGn\to
\widehat{PB_{n-1}}$ (see \cite{HK}). 
Namely, let 
$\widehat{B_{n-1}}=B_{n-1}\times_{PB_{n-1}}\widehat{PB_{n-1}}$. 
Then we can define a homomorphism $\Xi_n: \widehat{J_{n-1}}\to \widehat{B_{n-1}}$ 
by setting
$$
\Xi_n(\sigma_{p,q,r})=\beta_{[p,q],[q+1,r]}(\beta_{[p,p+r-q-1],[p+r-q,r]}\beta_{[p,q],[q+1,r]})^{-1/2}
$$ 
where $\beta_{[p,q],[q+1,r]}$ is the braid which interchanges the
intervals $[p,q]$ and \linebreak $[q+1,r]$. 


At the Lie algebra level, $\Xi_n$ defines a homomorphism of filtered Lie algebras
${\rm Lie}\hGn\to {\rm Lie}\widehat{PB_{n-1}}$.
Taking the associated graded of this map, 
and using the fact that ${\rm grLie}\hGn=\mathcal L_n\otimes\mathbb Q$,
we obtain the required statement.
\end{proof}

\begin{conjecture}\label{inje}
$\xi_n$ is injective.
\end{conjecture}

\begin{remarks} \begin{enumerate}
\item Conjecture \ref{inje} implies Conjecture
\ref{isomo1} (i) and Conjecture \ref{free}.
\item Conjecture \ref{inje} also implies Conjecture
\ref{freeness}, as we have natural morphisms $\mathbf L_{n+1}\to \mathbf
L_n$, whose kernels are known
to be free Lie algebras. On the other hand, by the Shirshov-Witt theorem,
a Lie subalgebra of a free Lie algebra is free.
\end{enumerate}
\end{remarks}


\section{Poset homology and a twisted version of $\Lambda_n$} \label{se:basis}

We now leave this world of Lie algebras and turn to the task of proving our main theorem describing the cohomology ring.  As a first step, we will examine our ring $ \Lambda_n $ more closely.  More precisely, we will consider a twisted version of $ \Lambda_n$ which has close connections to the homology of the poset of odd set partitions (which in turn has close connections with the operad of Lie 2-algebras \cite{HW}).


\subsection{Homology of the poset of odd set partitions}

Given a poset $P$ with bottom and top elements $\hat{0}$ and
$\hat{1}$, there is an associated chain complex
$\tilde C_*$ in which $\tilde C_{r+1}$ is the free $\ZZ$-module spanned by
chains
\[
(\hat{0}<x_1<\dots<x_r<\hat{1})
\]
with differential
\begin{align}
\partial_r&(\hat{0}<x_1<\dots<x_r<\hat{1})\notag\\
&=
\sum_{1\le i\le r} (-1)^{i-1}
(\hat{0}<x_1<\cdots<x_{i-1}<x_{i+1}<\cdots<x_r<\hat{1}).\notag
\end{align}
By our convention, if $\hat{0}=\hat{1}$ then $\tilde C_0(P)=\ZZ$ and $\tilde
C_r(P)=0$ for $r\ne 0$. The (shifted reduced) homology $\tilde H_*(P)$
is defined to be the homology of this chain complex.
If $P$ doesn't have a top element, we define $\tilde H_n(P):=\tilde H_{n+1}(P^+)$, where
$P^+=P\sqcup\{\hat{1}\}$.
Let $\tilde h_n(P)$ be the rank of $\tilde H_n(P)$.  
Note that we have shifted the degrees of the chains from the usual convention, in order for the following to hold.

\begin{prop}\label{Kun}
If $P$ and $Q$ are posets with bottom and top and
at least one of $\tilde H_*(P)$ and $\tilde H_*(Q)$ is free, then
\[
\tilde H_*(P\times Q) = \tilde H_*(P)\otimes \tilde H_*(Q)
\]
as graded modules.  In particular, the induced isomorphism $\tilde
H_*(P\times Q)\to \tilde H_*(Q\times P)$ is the standard exchange map for
the tensor product of supermodules.
\end{prop}

Here and below, as usual in homological algebra, the tensor product of supermodules is always taken
in the symmetric monoidal category of supermodules, which means that
the commutativity isomorphism is given by the formula $v\otimes w\to
(-1)^{\deg(v)\deg(w)}w\otimes v$ for homogeneous vectors $v,w$.

\begin{proof}
In terms of the order complex $|\Delta(P)|$ of $P$, we have
\[
\tilde H_*(P) = H_*(|\Delta(P)|,|\Delta(P)|\setminus\{x\}),
\]
where $x$ is the midpoint of the edge between $\hat{0}$ and $\hat{1}$ (the
local homology at that point).  The result then follows from the K\"unneth
formula for relative homology.
\end{proof}

Let $\Piodd{n}$ denote the poset of partitions of $\{1,2,\dots,n\}$ in
which each part has odd order.

\begin{theorem} \cite{Bj,CHR}
The poset $\Piodd{n}$ is totally semimodular, and thus Cohen-Macaulay;
that is, the homology of any closed interval in the poset is a free
$\ZZ$-module concentrated in the top degree.
\end{theorem}

In particular, $\tilde H_*(\Piodd{2n+1})$ is concentrated in degree
$n$, as is $\tilde H_*(\Piodd{2n+2})$; note that as the latter does not
naturally have a top element, it is thus implicitly added (and the degree shifted by 1).  
The further structure of these modules has been determined in \cite{CHR}.

\begin{theorem} \cite{CHR} \label{th:dimposethom}
We have the exponential generating functions
\[
\sum_{n\ge 0} \tilde h_n(\Piodd{2n+1}) \frac{u^{2n+1}}{(2n+1)!}
=
\arcsin(u)
\]
and
\[
\sum_{n\ge 0} \tilde h_n(\Piodd{2n+2}) \frac{u^{2n+2}}{(2n+2)!}
=
1-\sqrt{1-u^2}.
\]
\end{theorem}


In \cite{Koz}, Kozlov considers the following spectral sequence.  Let $P$
be a Cohen-Macaulay poset.
The rank of an element $x\in P$ is then the unique $r$ such that $\tilde H_r([\hat 0,x])\not =0$
(we take by
convention that $\rank(\hat{0})=0$).
Consider the filtration of the associated chain complex given by
\[
F^rC_{k+1} = \langle (\hat{0}<x_1<\dots<x_k<\hat{1}):\rank(x_k)\le r\rangle.
\]
Namely, for our purposes, we consider chains of set partitions such that the
coarsest (=biggest) partition in the chain has at least $n-2r$ parts.  Then Kozlov
observes that the associated spectral sequence looks like
\[
E^1_{k,l} = \bigoplus_{\rank(x)=l} \tilde H_k([\hat{0},x])\,\,\Rightarrow\,\,\tilde H_{k+1}(P),
\]
where the indices $k,l$ run from 0 to $\rank(\hat 1)-1$.
But this is 0 unless $k=l$, and thus the spectral sequence collapses on the
$E^2$ page.  Since the spectral sequence converges to $\tilde H_*(P)$ and
this is concentrated in the top degree, we find that $E^2_{k,l}=0$ unless
$k=l=\rank(\hat 1)-1$, in which case it equals the homology of $P$.

In other words, we obtain the following result.

\begin{thm}\label{thm34}
If $P$ is a Cohen-Macaulay poset, the Whitney homology groups
\[
W_k(P) = \bigoplus_{x\in P} \tilde H_k([\hat{0},x])
       = \bigoplus_{\rank(x)=k} \tilde H_*([\hat{0},x])
\]
form a canonical exact sequence
\[
0\to W_{\rank(\hat{1})}\to\cdots\to W_1\to W_0\to 0.
\]
\end{thm}

\begin{remark}
Note that $W_{\rank(\hat{1})}=\tilde H_*(P)$ and $W_0=\ZZ$.  If $P$ has a
bottom but not a top, we by convention do not include the $\tilde H_*(P)$
term in $W_*(P)$. 
\end{remark}

By construction, this exact sequence respects products when the reduced
homology modules are free.  The differential corresponding to
this chain complex acts on saturated chains (in which each step is a
covering relation) by
\[
\partial_W(\hat{0}<x_1<\dots<x_r)
=
(-1)^{r-1}
(\hat{0}<x_1<\dots<x_{r-1}).
\]
Note that since $P$ is Cohen-Macaulay, only saturated chains appear in the Whitney
homology.

\subsection{A twisted version of $\Lambda$}

For a finite set $S$, let $\tilde\Lambda(S)$ be the supercommutative
$\ZZ$-algebra with symmetric degree 1 generators $\tau_{ijk}$ for $i$, $j$,
$k$ distinct elements of $S$, subject to the quadratic relations
\begin{gather}
\tau_{ijk}\tau_{jkl} = 0\notag\\
\tau_{ijk}\tau_{klm}+ \tau_{jkl}\tau_{lmi}+ \tau_{klm}\tau_{mij}+
\tau_{lmi}\tau_{ijk}+ \tau_{mij}\tau_{jkl}=0.\notag
\end{gather}

This is of course quite closely related to the $S_n$-invariant presentation
of $\Lambda_{n+1}$ considered above, except for the change in symmetry of
the generators.  In fact, experimentation for small $n$ immediately
suggests that the two algebras have the same Hilbert series, a fact which
we will confirm below.

A triangle graph is pair $ (V, E) $ where $ V $ is a finite set and $ E $ a subset of the set of three element subsets of $ V $.  $ V $ is called the set of vertices and $ E $ the set of edges.  More generally, there is the notion of hypergraph, where subsets of any size are allowed.  

Many notions from graph theory may be extended to triangle graphs.  A path of length $ n $ in a triangle graph, from vertex $ x $ to a vertex $ y $, is a sequence of edges $ e_1, \dots, e_n $ and vertices $ x= v_0, \dots, v_n=y $ with $ v_i \in e_i \cap e_{i+1} $ for $ 0 < i < n $, $ v_0 \in e_1 $ and $ v_n \in e_n$.  A cycle is a path of length greater than 1 from a vertex to itself which consists of distinct edges and vertices.  A triangle forest is a triangle graph without any cycles and a connected triangle forest is called a triangle tree.

Given a monomial in the triples $\tau_{ijk}$, there is a natural triangle graph with vertex set $S $ and with an edge for each variable.  We can consider the associated partition of the set $S$ into connected components.   This gives a grading of the exterior algebra with respect to the lattice of partitions.  That is, the product of the homogeneous components associated to two partitions $\pi_1$, $\pi_2$ lies in the homogeneous component associated to their join.  Since the ideal of relations of
$\tilde\Lambda(S)$ is homogeneous with respect to this grading, $\tilde\Lambda(S)$ inherits this
grading. Thus for each partition $\pi$ of $S$, we associate a
corresponding space $\tilde\Lambda[\pi]$ and we have
$\tilde\Lambda(S)=\underset{\pi}{\bigoplus}\,\tilde\Lambda[\pi]$.  

\begin{gather}\label{AnExample}
\begin{matrix}
\psfig{file=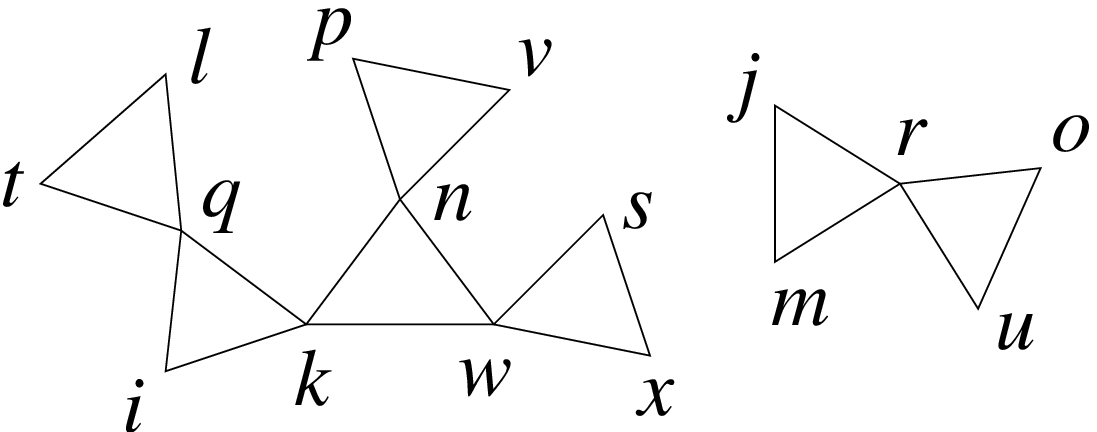,width=5cm}
\end{matrix}
\\
\nonumber
\text{\it The triangle graph for the monomial $\tau_{tlq}\tau_{pvn}\tau_{qik}\tau_{knw}\tau_{wsx}\tau_{mjr}\tau_{rou}$.}
\end{gather}

\begin{lemma}
If $\pi_1$ and $\pi_2$ are partitions of disjoint sets $S_1$ and $S_2$,
then
\[
\tilde\Lambda[\pi_1\cup\pi_2] = \tilde\Lambda[\pi_1]\bigotimes
\tilde\Lambda[\pi_2]
\]
under the natural imbedding of
$\tilde\Lambda(S_1)\bigotimes\tilde\Lambda(S_2)$ in
$\tilde\Lambda(S_1\cup S_2)$.
\end{lemma}

In particular, the structure of these homogeneous components is determined
by the structure for partitions with a single part. 

\subsection{Spanning the algebra $ \tilde\Lambda(S)$}
Direct calculation
gives some small cases of $\tilde\Lambda[\pi]$:
\begin{align}
\tilde\Lambda[\{1\}]&=\langle 1\rangle\notag\\
\tilde\Lambda[\{1,2\}]&=0\notag\\
\tilde\Lambda[\{1,2,3\}]&=\langle \tau_{123}\rangle\notag\\
\tilde\Lambda[\{1,2,3,4\}]&=0.\notag
\end{align}
For $n=5$, we find the following.

\begin{lemma}
The space $\tilde\Lambda[\{1,2,3,4,5\}]$ is spanned by the nine products
$\tau_{12i}\tau_{jkl}$ for which $|\{1,2,i,j,k,l\}|=5$.
\end{lemma}

\begin{lemma}\label{lem:path_basic}
For any pair of distinct elements $j,j'\in S$,
$\tilde\Lambda[S]$ is spanned by monomials in which at least one triple
contains $\{j,j'\}$.
\end{lemma}

\begin{proof}
Any monomial in $\tilde\Lambda[S]$ corresponds to a connected
triangle graph, and thus that triangle graph will contain a shortest path from
$j$ to $j'$.  We proceed by induction on the length of that path, noting
that if the path has length 1, we are done.

Assuming therefore that the length is greater than 1, let the shortest path
begin by $j$, $k$, $l$.  Then the monomial contains triples
\[
\tau_{ijk}\tau_{klm}
\]
for some $i,m\in S$.
If we apply the quadratic relation, the resulting monomials all contain a
triple joining $j$ and $l$, and thus the length of the shortest path
from $j$ to $j'$ in the new monomials has been decreased by 1.
\end{proof}

\begin{corollary}
If the triangle graph of a monomial in $\tilde\Lambda[S]$ contains a cycle,
then that monomial is 0.
\end{corollary}

\begin{proof}
The same argument shows that we can reduce the cycle to length 2.  But a
cycle of length 2 is of the form
$\tau_{ijk}\tau_{ijl}=0$.
\end{proof}

In other words, the triangle graph of any nonzero monomial in $\tilde\Lambda(S)$ is a forest.
This immediately implies that the grading by partitions refines the grading
by degree.

\begin{theorem}
Let $\pi$ be a partition of a set $S$ of cardinality $n$. Then any element of
$\tilde\Lambda[\pi]$ has degree $(n-|\pi|)/2$.  In particular, if $n-|\pi|$
is odd, then $\tilde\Lambda[\pi]=0$.
\end{theorem}

\begin{proof}
Consider a nonvanishing monomial of degree $d$ in $\tilde\Lambda[\pi]$.
The incidence graph of its associated triangle graph has $n+d$ vertices and
$3d$ edges; since it is a forest, it must therefore have $n-2d$ connected
components.
\end{proof}

\begin{corollary}
If $\pi\not\in\Piodd{n}$, then $\tilde\Lambda[\pi]=0$.
\end{corollary}

\begin{proof}
Indeed, the components of a triangle forest all have odd order.
\end{proof}

Standard results on exponential generating functions then imply the following.

\begin{corollary}\label{cor311}
Let $f(u,t)$ be the generating function
\[
\sum_{0\le k}
\dim(\tilde\Lambda[\{1,2,\dots,2k+1\}])
t^k \frac{u^{2k+1}}{(2k+1)!},
\]
and let $P(u,t)$ be the generating function
\[
\sum_{0\le n} P_n(t) \frac{u^n}{n!},
\]
where $P_n(t)$ is the Poincar\'e polynomial of $\tilde\Lambda_n$. 
Then $P(u,t) = \exp(f(u,t))$.
\end{corollary}

Using Lemma \ref{lem:path_basic} to its fullest, we can obtain the
following spanning set of $\tilde\Lambda(S)$.  Choose a total ordering on
$S$. 

\begin{definition}
We define {\bf basic} triangle trees and forests as follows, by
induction on the number of triangles.

\begin{itemize}
\item[1.] A single point is a basic triangle tree.
\item[2.] A nontrivial triangle tree is basic iff the two smallest
  vertices are on a common triangle, and each of the three components after
  removing that triangle is basic.
\item[3.] A triangle forest is basic iff each component is basic.
\end{itemize}
\end{definition}

For example, (\ref{AnExample}) is a basic triangle forest with respect to the usual order on the letter $i,j,\ldots,x$.

\begin{proposition}\label{prop312}
The algebra $\tilde\Lambda(S)$ is spanned by the monomials associated to
basic triangle forests.
\end{proposition}

\begin{proof}
This is a simple induction; by Lemma \ref{lem:path_basic}, any monomial associated to a tree can
be expanded in terms of monomial associated to trees in which the two smallest vertices are on a
common triangle, and similarly for the monomials associated to the
subtrees.
\end{proof}

If $h(u)$ is the exponential generating function for basic trees, then we
claim that $h''(u)=u h'(u)^3$.  Indeed, the left-hand side is the
exponential generating function for basic trees on $n+2$ vertices, while
the right-hand side corresponds to the following:
\begin{itemize}
\item[a.] A choice of element $k\in \{3,\dots,n+2\}$
\item[b.] An ordered set partition of $\{3,\dots,n+2\}\setminus \{k\}$ into
  three parts $P_1$, $P_2$, $P_3$.
\item[c.] Basic trees on $P_1\cup\{1\}$, $P_2\cup\{2\}$, $P_3\cup\{k\}$.
\end{itemize}
Adjoining the triangle $\{1,2,k\}$ to the three basic trees gives a
basic tree on $n+2$ vertices, establishing the differential equation.
The unique solution with $h(0)=0$, $h'(0)=1$ is $h(u)=\arcsin(u)$.

We will now show that basic triangle trees form a basis of $\tilde\Lambda[S]$ 
(and similarly for basic triangle forests and $\tilde\Lambda(S)$).
As a consequence, $\arcsin(u)$ is also the exponential generating function for $\dim\tilde\Lambda[S]$.

\subsection{A basis for $ \tilde\Lambda(S)$}
\begin{thm}\label{thm:homology_algebra}
For a partition $\pi\in\Piodd{n}$, let $S_\pi\subset S_n$ be the corresponding product of symmetric groups.
Then there is a canonical
isomorphism of graded $\ZZ[S_\pi]$-modules between $\tilde\Lambda[\pi]$ and
$\tilde H_*([\hat{0},\pi])$.

Moreover, these abelian groups have bases indexed by basic triangle forests with component partition $\pi$.
\end{thm}

\begin{proof}
Since both $\tilde\Lambda[\pi]$ and $\tilde H_*([\hat{0},\pi])$ are
multiplicative under disjoint union of partitions, it is enough to treat
the case $|\pi|=1$, namely $[\hat 0,\pi]=\Piodd{n}$.

Given a triangle tree, any ordering of the triangles gives rise
to a chain of set partitions; the $k$-th partition in the chain is
the set of components of the subgraph spanned by the first $k$
triangles.  We thus obtain a map $\tilde\Lambda[\hat 1] \to
\tilde C_*(\Piodd{n})$, taking a tree to the alternating
sum over the chains resulting from all orderings of the triangles.
The five-term relation (and any tree multiple thereof) is
annihilated by this map, so it is indeed well-defined.

We claim that the image of a tree is a cycle, and we thus obtain a
well-defined map to $\tilde H_*(\Piodd{n})$.  Indeed, the operation of
removing the $l$-th partition in a chain is nearly invariant under
swapping the $l$-th and $l+1$-st triangles, except that swapping the
triangles gives rise to an overall sign.

Since the number of basic trees is equal
to the rank of the free $\ZZ$-module $\tilde H_*(\Piodd{n})$ (by the remark following Proposition \ref{prop312} and
Theorem \ref{th:dimposethom}), the result will follow if we can exhibit a
set of cochains such that the induced pairing with basic trees is
triangular with unit diagonal.  Note from the description of the map from
trees to chains that the pairing of an elementary cochain (corresponding
to a maximal chain in $\Piodd{n}$) with a triangle tree is $\pm 1$ or $0$.
Moreover, the pairing is nonzero precisely if there exists an ordering of
the triangles in the tree corresponding to the same chain of odd set
partitions.

First, some additional terminology.  The ``root'' of a basic tree is the
triangle containing the two minimal elements of its support.  If we remove
the root, we obtain three components; the one {\em not} containing one of
the original minimal elements will be called the ``stepchild''.

Given a basic tree $T$, $\text{size}(T)$ represents the size of its support,
and we define
\[
\rank(T):=\frac{\text{size}(T)-\text{size}(\text{stepchild}(T))}{2}.
\]
Then if $\text{size}(T)=2n+1$, we can associate a composition $\mu(T)$ of $n$
with first part $\rank(T)$, second part the rank of its stepchild, and so
forth.  This then induces a partial order on basic trees by
lexicographically ordering the associated compositions.

Finally, we define the ``keystone'' of a basic tree as follows.  If the
stepchild of $T$ is a single vertex, its keystone is simply the root.
Otherwise, the keystone of $T$ is the keystone of its stepchild.

\begin{lemma}\label{lem314}
Let $T$ be a basic triangle tree and $F$ be the forest 
obtained by removing its keystone.  
Then any other basic tree $T'$ containing $F$ has $\mu(T')>\mu(T)$.
\end{lemma}

\begin{proof}
If the root of $T$ was the keystone, so $T$ had maximal composition $n$,
then $F$ consists of the component of 1, the component of 2, and a single
isolated vertex $k$, and thus the only basic tree containing $F$ is obtained by
adjoining the triangle $12k$, so the result holds in this case.

Otherwise, the root of $T$ is still in $F$.  If we remove it, we obtain a five
tree forest, consisting of the component of 1, the component of 2, and
three other components all contained in the stepchild of $T$.  
Now, $T'$ is obtained from $T$ by removing the keystone and adding another triangle instead.
If that triangle meets the component of 1 or 2, the composition clearly
increases.  Otherwise, the result follows by induction applied to the
stepchild of $T$.
\end{proof}

Now, for each basic forest $F$, we choose a triangle $t_F$ which is the
keystone of some component.  This then inductively gives rise to an
ordering $(t_i)_{i=1}^n$ on the triangles of any tree $T$ by letting
$t_n:=t_T$ and $t_{i-1}:=t_{T\setminus\{t_i\ldots t_n\}}$.  Let $\alpha_T$
be the sequence of set partitions obtained by successively adding $t_1,
t_2,$ etc, which we identify with the corresponding elementary cochain.

By Lemma \ref{lem314}, we can reconstruct the sequence of triangles from
$\alpha_T$ by taking at each step the triangle that minimizes the
composition associated to the resulting tree.  In other words, $T$ is
reconstructed from $\alpha_T$ as the minimal tree having a nonzero (and
thus unit) pairing with $\alpha_T$. This immediately gives us the desired
triangularity.
\end{proof}

\begin{corollary}\label{cor:hilser_t}
Basic forests form a basis of $\tilde\Lambda(S)$.
The corresponding exponential generating function (see Corollary \ref{cor311}) is given by
\[
P(u,t)=\exp\big(\arcsin(u\sqrt{t})/\sqrt{t}\big),
\]
and the Hilbert series of $\tilde\Lambda(\{1,2,\dots,n\})$ is
\[
P_{n+1}(t)=\prod_{0\le k<(n-2)/2}(1+(n-2-2k)^2t).
\]
\end{corollary}

\begin{proof}
The only thing remaining to show is that the exponential generating
function implies the Hilbert series, or in other words that
\[
\exp\big(\arcsin(u\sqrt{t})/\sqrt{t}\big)
=
\sum_{n\ge 0}\:
\prod_{0\le k<(n-2)/2} (1+(n-2-2k)^2t)
\frac{u^n}{n!}.
\]
To this end, we observe that if $t=-1/k^2$ for a positive integer $k$, then
\begin{equation}\label{eq:p22}
\exp\big(\arcsin(u\sqrt{t})/\sqrt{t}\big)
=
\big(u/k + \sqrt{1+u^2/k^2}\big)^k.
\end{equation}
But (\ref{eq:p22}) is the sum of a polynomial of degree $k$ and of a function satisfying $f(u)=(-1)^{k+1}f(-u)$.
Thus, if $n-k$ is an even positive integer, we have $\frac{d^n}{du^n}P(u,-1/k^2)|_{u=0}=0$.
The coefficient of $u^n/n!$ in $P(u,t)$ is therefore divisible by $(1+(n-2i)^2t)$ for $1\le i< n/2$.
Since we know the degree of $P_{n+1}(t)$,
this determines the above expansion up to a constant which can be set via
the limit $P(u,0)=\exp(u)$.
\end{proof}

\begin{cor}
There is a canonical isomorphism
\begin{equation}\label{abi}
\tilde\Lambda(\{1,2,\dots,n\})\cong W_*(\Piodd{n}).
\end{equation}
\end{cor}

It will be of use to know how the canonical differential on $W_*$
is expressed on $\tilde\Lambda(\{1,2,\dots,n\})$.

\begin{lem}\label{lem317}
Under the isomorphism (\ref{abi}), the canonical differential on $W_*$ takes
a forest monomial
\[
\tau_{a_1b_1c_1}\tau_{a_2b_2c_2}\dots\tau_{a_lb_lc_l}
\]
to
\[
\partial(\tau_{a_1b_1c_1}\tau_{a_2b_2c_2}\dots\tau_{a_lb_lc_l})
=
\sum_{1\le i\le l}
(-1)^{i-1} \prod_{j\ne i} \tau_{a_jb_jc_j}.
\]
\end{lem}

\begin{proof}
The action of the canonical differential is to remove the last
element in the chain of partitions.  This is equivalent to removing the
last triangle in the ordering.  Summing over the possible last triangles
gives the stated result.
\end{proof}

This is the differential for the reduced homology of a simplicial complex, with
a vertex for each triangle and a simplex for each basic forest.

\begin{cor}
The simplicial complex of basic forests on $\{1,2,\dots,2n+1\}$ is
contractible.  The homology of the simplicial complex of disconnected
basic forests on $\{1,2,\dots, n\}$ is isomorphic to the homology of
$\Piodd{n}$.
\end{cor}

\section{The cohomology of the Bockstein} \label{se:bound}
As a key step in determining the cohomology ring of $ M_n$, we will compute
the cohomology of $ M_n$ with coefficients in $ \ZZ/4\ZZ $, modulo 2-torsion.
To do so we will start with the cohomology ring of the complex moduli space
$M_n^\CC:=\overline{M_{0,n}}(\CC)$.  In Theorem \ref{thm45}, we will use
Keel's work to find an explicit basis for $H^*(M_n^\CC, \ZZ) $ which is
convenient for our purposes.  We use this cohomology ring to determine
the cohomology of the real moduli space $ M_n $ with coefficients in $
\FF_2$ (Theorem \ref{th:cohommod2}).  This mod 2 cohomology $ H^*(M_n,
\FF_2) $ has a differential, the Bockstein, whose cohomology is the
cohomology of $ M_n $ with coefficients in $ \ZZ/4\ZZ$, modulo 2-torsion.
In Theorem \ref{thm52}, we compute the cohomology of this Bockstein
differential and show that it is isomorphic to $ \tilde\Lambda(\{1, \dots,
n-1\}) \otimes \FF_2$.  Thus the resulting Betti numbers give us the
desired upper bound on the Betti numbers of $ M_n $ (Corollary
\ref{th:bound}).

\subsection{Cohomology of $ M_n^\CC $}
We begin by recalling the following presentation, due to Keel, of
$H^*(M_n^\CC,\ZZ)$.  For our purposes, it will be useful to single out one
of the marked points, and thus consider $M_{n+1}^\CC$, with points marked
$0$, $1$, \dots, $n$.

\begin{theorem}\cite{Keel}
The commutative ring $H^*(M_{n+1}^\CC,\ZZ)$ is generated by elements (of
degree 2) $D_S$, one for each subset $S\subset \{0,1,2,\dots,n\}$ with $2\le
|S|\le n-1$, subject to the following relations.
\begin{itemize}
\item[1.] $D_S = D_{\{0,1,\dots,n\}\setminus S}$.
\item[2.] For distinct elements $i,j,k,l\in \{0,1,\dots,n\}$,
\[
\sum_{\substack{i,j\in S\\k,l\notin S}} D_S =
\sum_{\substack{i,k\in S\\j,l\notin S}} D_S.
\]
\item[3.] If $S\cap T\notin \{\emptyset,S,T\}$ and $S\cup
  T\ne\{0,1,\dots,n\}$, then $D_S D_T=0$.
\end{itemize}
The exponential generating function
\[
A(u,t) := \sum_{n\ge 2}\, \sum_{k=0}^{n-2} \dim(H^{2k}(M_{n+1}^\CC,\ZZ))
t^k \frac{u^n}{n!}
\]
satisfies the differential equation
\begin{equation}\label{bti}
\frac{\partial}{\partial u} A(u,t) = u + (1+t)A(u,t) + t A(u,t) \frac{\partial A(u,t)}{\partial u}.
\end{equation}
\end{theorem}

\begin{remark}
The class $D_S$ has a geometrical interpretation as the class of the
divisor of $M_{n+1}^\CC$ consisting of singular genus 0 curves in which
the removal of a singular point separates the points in $S$ from the points
not in $S$.
\end{remark}

For our purposes, we will need an alternate presentation involving only $S\subset\{1,\ldots,n\}$.
First, define an element
\[
D_{\{1,2,\dots,n\}}
=
-\sum_{\{1,2\}\subset S\subsetneq \{1,2,\dots,n\}} D_S.
\]

\begin{lem}
The element $D_{\{1,2,\dots,n\}}$ is invariant under the action of $S_n$.
\end{lem}

\begin{proof}
We need simply to show that
\[
\sum_{\{1,2\}\subset S\subsetneq \{1,2,\dots,n\}} D_S
=
\sum_{\{1,3\}\subset S\subsetneq \{1,2,\dots,n\}} D_S.
\]
But if we eliminate the common terms from both sides, this becomes
\[
\sum_{\substack{1,2\in S\\0,3\notin S}} D_S
=
\sum_{\substack{1,3\in S\\0,2\notin S}} D_S,
\]
which holds by the linear relation.
\end{proof}

\begin{remark}
It follows from Proposition 1.6.3 of \cite{KM} that
\[
D_{\{1,2,\dots,n\}} = -c_1(\ell_0),
\]
where $\ell_0$ is the line bundle obtained by taking the tangent line at
the $0$-th marked point, and $c_1$ is the first Chern class.
%
\end{remark}

Now, for $2\le |S|\le n$, define
\[
\Pi_S = -\sum_{S\subset T\subset\{1,2,\dots,n\}} D_T.
\]
Since
\[
D_S = -\sum_{S\subset T\subset \{1,2,\dots,n\}} (-1)^{|T|-|S|}\Pi_T,
\]
these elements span $H^*(M_{n+1}^\CC,\ZZ)$.  On the other hand, by the
definition of $D_{\{1,2,\dots,n\}}$, we find that $\Pi_S=0$ whenever $|S|=2$,
and thus $H^1(M_{n+1}^\CC,\ZZ)$ is spanned by the elements $\Pi_S$ for
$|S|\ge 3$.  There are $2^n-\binom{n}{2}-n-1$ such elements, which equals
the rank of $H^1(M_{n+1}^\CC,\ZZ)$, and thus the elements $\Pi_S$ for
$|S|\ge 3$ form a basis of $H^1(M_{n+1}^\CC,\ZZ)$, and generate the
cohomology ring.

\begin{proposition}
The elements $\Pi_S$ satisfy the following quadratic relation.  For any
subsets $S$, $T\subset\{1,2,\dots,n\}$ such that $S\cap
T\notin\{\emptyset,S,T\}$,
\begin{equation}\label{hrt}
(\Pi_S-\Pi_{S\cup T})(\Pi_T-\Pi_{S\cup T})
=
0.
\end{equation}
\end{proposition}

\begin{proof}
Note that the conditions on $S$ and $T$ imply $2\le |S|,|T|\le n-1$.
We have
\[
\Pi_S-\Pi_{S\cup T} = \sum_{\substack{S\subset U\\T\not\subset U}} D_U\quad\qquad\text{and}\quad\qquad
\Pi_T-\Pi_{S\cup T} = \sum_{\substack{T\subset V\\S\not\subset V}} D_V.
\]
In particular, for any terms $D_U$ and $D_V$ in the respective sums, we
have $S\cap T\subset U\cap V$, so $U\cap V\ne \emptyset$, and
$S,T\not\subset U\cap V$, so $U\cap V\notin \{U,V\}$.  In other words,
$D_UD_V=0$ for any such pair, and the product of the two sums vanishes
termwise.
\end{proof}

\begin{remark}
When $|T|=2$, we have $\Pi_T=0$. Since $S\cup T=S\cup\{i\}$ for some $i\notin S$,
\begin{equation}\label{bw}
\Pi_{S\cup \{i\}}(\Pi_S-\Pi_{S\cup \{i\}})=0
\end{equation}
is a special case of (\ref{hrt}).
More generally, it follows from an easy induction that for any disjoint
sets $S$, $T$,
\begin{equation}\label{bbw}
\Pi_{S\cup T}^{|S|+1} = \Pi_{S\cup T}^{|S|}\Pi_T.
\end{equation}
\end{remark}

We claim that (\ref{hrt}) are the only relations satisfied by the
$\Pi_S$.  To prove this, it will in fact be simplest to give a basis for the
ring thus presented, and show that it is a free $\ZZ$-module with the
correct Hilbert series.

In fact, we can give a Gr\"obner basis for the ring, which will moreover be
invariant under the $S_n$ action.  We need one more set of relations, which
can be deduced from (\ref{hrt}).

\begin{lem}
The elements $\Pi_S$ satisfy the following relation for all $k\ge 0$.  Let
$S_0$, $S_1$, \dots, $S_k$ be disjoint sets, with union $S$; suppose
moreover that $|S_i|\ge 3$ for $1\le i\le k$.  Then
\begin{equation}\label{dxd}
\Pi_S^{|S_0|+k-1} \prod_{1\le i\le k} (\Pi_{S_i}-\Pi_S)
=
0.
\end{equation}
\end{lem}

\begin{proof}
For $k=0$, this becomes the statement
\begin{equation}\label{sot}
\Pi_S^{|S|-1}=0,
\end{equation}
which follows from (\ref{bw}) and the fact that $\Pi_T=0$
if $|T|=2$.

For $k=1$ the claim is equivalent to (\ref{bbw}), so let us assume that $k\ge 2$.
By applying (\ref{bbw}) multiple times, we get
\[
\Pi_S^{|S_0|+k-1} \prod_{1\le i\le k} (\Pi_{S_i}-\Pi_S)
=
\Pi_S^{|S_0|}
\Pi_{S\setminus S_0}^{k-1} \prod_{1\le i\le k} (\Pi_{S_i}-\Pi_{S\setminus S_0})
\]
and may thus assume $S_0=\emptyset$.

For $k=2$, choose $i\in S_2$ and consider the known relation
\begin{equation}\label{dcr}
(\Pi_{S_1}-\Pi_{S_1\cup \{i\}})\Pi_{S_2}\Pi_{S_1\cup \{i\}}=0.
\end{equation}
By (\ref{hrt}), we have
\[
\Pi_{S_2}\Pi_{S_1\cup \{i\}}
=
\Pi_{S_2}\Pi_{S_1\cup S_2}
+
\Pi_{S_1\cup \{i\}}\Pi_{S_1\cup S_2}
-
\Pi_{S_1\cup S_2}^2,
\]
so we can simplify (\ref{dcr}) to
\[
\begin{split}
0&=(\Pi_{S_1}-\Pi_{S_1\cup\{i\}})(\Pi_{S_2}-\Pi_{S_1\cup S_2})\Pi_{S_1\cup S_2}\\
&=
(\Pi_{S_1}-\Pi_{S_1\cup S_2})(\Pi_{S_2}-\Pi_{S_1\cup S_2})\Pi_{S_1\cup S_2}
\end{split}
\]
as required.

For $k>2$, let $T=S_1\cup S_2$.  We need to show that
\begin{equation}\label{On}
\Pi_S^{k-1} \prod_{1\le i\le k} (\Pi_{S_i}-\Pi_S)
=
0.
\end{equation}
If we subtract from (\ref{On}) the relation
\[
(\Pi_{S_1}-\Pi_T)(\Pi_{S_2}-\Pi_T)\Pi_T
\Pi_S^{k-2}\prod_{i>2} (\Pi_{S_i}-\Pi_S),
\]
we find by setting $\Pi_T=\Pi_S$ that the result is divisible by $\Pi_T-\Pi_S$.
It is therefore divisible by
\[
\Pi_S^{k-2}(\Pi_T-\Pi_S)\prod_{i>2} (\Pi_{S_i}-\Pi_S),
\]
and so is 0 by induction.
\end{proof}

\begin{theorem}\label{thm45}
The algebra $H^*(M_{n+1}^\CC,\ZZ)$ is freely spanned by monomials of the form
$\prod_{|S|\ge 3} \Pi_S^{d_S}$ satisfying the following conditions.
\begin{itemize}
\item[1.] If $d_S>0$, $d_T>0$, then $S\cap T\in \{\emptyset,S,T\}$.
\item[2.] For each $S$ such that $d_S>0$, let $S_1$,\dots, $S_k$ be the
  maximal proper subsets of $S$ such that $d_{S_i}>0$, disjoint by
  condition 1.  Then
\[
d_S<k-1+|S|-\sum_i |S_i|.
\]
\end{itemize}
Equivalently, the relations (\ref{hrt}) and (\ref{dxd}) form a Gr\"obner
basis of (the ideal of relations of) the cohomology ring 
$H^*(M_{n+1}^\CC,\ZZ)$ with
respect to the grevlex order, relative to any ordering on the variables
$\Pi_S$ extending inclusion.
\end{theorem}

\begin{proof}
Any monomial not satisfying the above conditions can be expanded in smaller
monomials using the given relations; as a result, the ``good'' monomials
span.  The theorem will thus follow if we can show that they form a basis.

For this, we need simply count the monomials and show that we obtain the
correct exponential generating function.  Now, consider the exponential
generating function $C(u,t)$ of monomials with $d_{\{1,2,\dots,n\}}>0$,
extended to include the constant $1$ for $n=1$.  By standard manipulations
of exponential generating functions, this satisfies
\[
C(u,t) = u + \sum_{m\ge 3}\, \sum_{l=1}^{m-2} t^l C(u,t)^m/m!,
\]
where $m=k+|S|-\sum_i |S_i|$ and $l=d_S$.  If we multiply by $1-t$
and simplify, this becomes the functional equation
\[
\exp(tC(u,t)) = 1+tu+t^2(\exp(C(u,t))-1-u).
\]

Now, the exponential generating function of all monomials (including $1$
for $n\le 2$) is $B(u,t):=\exp(C(u,t))$, and thus satisfies
\[
B(u,t)^t = 1+tu+t^2(B(u,t)-1-u).
\]
Differentiating with respect to $u$, we get
\[
\frac{\partial B(u,t)}{\partial u} = \frac{B(u,t)}{1-t(B(u,t)-1-u)}.
\]
Therefore $B(u,t)-1-u$ satisfies the differential equation (\ref{bti}),
and we indeed have the correct number of monomials.
\end{proof}

\begin{remark}
A similar $\ZZ$-basis of $H^*(M_{n+1}^\CC,\ZZ)$ was given by Yuzvinsky
\cite{Yuz}, based on the generators $D_S$ for $S\subset\{1,2,\dots,n\}$,
$|S|>2$, but otherwise the same.  It differs in two important respects:
first, while it respects a filtration by set partitions, it does not
respect the grading by set partitions; and second, the associated Gr\"obner
basis is more complicated.  On the other hand, Yuzvinsky's basis is more
amenable to computation of the Poincar\'e pairing.
\end{remark}

One important consequence of this presentation is that there is a natural
grading of $H^*(M_{n+1}^\CC,\ZZ)$ indexed by the lattice of set
partitions of $\{1,2,\dots,n\}$.  
Indeed, the hypergraphs associated to each monomial of each
relation all have the same set of connected components; that set of
connected components is the associated set partition.
By inspection of the canonical monomials, we also find that the partitions
that appear never have sets of size 2.

\subsection{Cohomology with coefficients in $ \FF_2$}

The above presentation for $H^*(M_{n+1}^\CC,\ZZ)$ immediately gives rise
to a presentation for $H^*(M_{n+1}^\RR,\FF_2)$, by the following result.

\begin{theorem} \label{th:cohommod2}
There is a canonical isomorphism between $H^{2k}(M_{n+1}^\CC,\FF_2)$
and $H^k(M_{n+1}^\RR,\FF_2)$ for each $k$, compatible with the ring
structures.
\end{theorem}

\begin{proof}
Indeed, Keel \cite{Keel} showed that
$H^*(M_{n+1}^\CC,\ZZ)$ is generated by algebraic divisors, all of which
are in fact rational over $\RR$.  It follows that $\overline{M_{0,n+1}}$ is an
``algebraically maximal'' variety in the sense of \cite{Kr} and thus the
isomorphism follows from the main result of that paper.
\end{proof}

\begin{remark}
Januszkiewicz (personal communication) has suggested another proof.
One can construct $M_{n+1}$ as an iterated blowup of projective $n$-space
along the standard $A_{n-1}$ hyperplane arrangement.  The desired agreement
of mod 2 cohomology holds for $\RP^n$ and $\CP^n$ and is preserved under
blowing up of a linear subspace, and the result follows by induction.

See also \cite{HHP} for a more general class of spaces satisfying this same relation to their real locus.
\end{remark}

We denote the piece of $H^k(M_{n+1}^\RR,\FF_2)$ indexed by the set
partition $\pi$ by $\Xi^k[\pi]$; if $\pi$ is a partition of a subset of
$\{1,2,\dots,n\}$, we adjoin singletons as necessary to make it a
partition.  We will also use $\Xi^*_n$ to denote the entire cohomology
ring $H^*(M_{n+1}^\RR,\FF_2)$.

\subsection{The Bockstein map and its cohomology}
For any space $ X $, consider the Bockstein map $\beta : H^*(X, \ZZ/2\ZZ) \rightarrow H^*(X, \ZZ/2\ZZ) $, which is the connecting map for the 
exact sequence
\[
0\to \ZZ/2\ZZ\to \ZZ/4\ZZ\to \ZZ/2\ZZ\to 0.
\]
The Bockstein is a differential. From the long exact sequence in 
cohomology, we see that its cohomology is isomorphic to the quotient of the mod 4 cohomology by the 2-torsion,
\begin{equation*}
H_\beta^*(X) \cong H^*(X, \ZZ/4\ZZ)/H^*(X, \ZZ/4\ZZ)\langle 2 \rangle.
\end{equation*}

On goal is to compute the cohomology of the Bockstein map acting on $H^*(M_{n+1}, \ZZ/2\ZZ) = \Xi^*_n$ in order to determine $$ H^*(M_{n+1}, \ZZ/4\ZZ)/H^*(M_{n+1}, \ZZ/4\ZZ)\langle 2 \rangle. $$  We denote this cohomology $ H^*_\beta(\Xi^*_n) $.  The action of $\beta$ can be determined from the
following two properties (corresponding to the fact \cite[Sec.4.L]{Hat} that $\beta$
is the first Steenrod square).
\begin{itemize}
\item[1.] $\beta$ is a derivation of degree 1.
\item[2.] If $x\in \Xi^1_n$, then $\beta(x)=x^2$.
\end{itemize}
Since $\Xi^*_n$ is generated in degree 1, this uniquely determines the
action of $\beta$ on the entire algebra:
\[
\beta(v) = \sum_S \Pi_S^2 \frac{\partial}{\partial \Pi_S} v.
\]
(This is well-defined since $\beta(I)\subset I$ where $I$ is the defining
ideal; in fact, if $r$ is one of the defining relations (\ref{hrt}), then $\beta(r)$ is
a multiple of $r$.)  In particular, $\beta$ is homogeneous with respect to
the grading by set partitions, and thus
\[
H^*_\beta(\Xi^*_n) = \bigoplus_\pi H^*_\beta(\Xi^*[\pi]).
\]
Moreover,
\[
\Xi^*[S_1,S_2,\dots,S_k] = \bigotimes_i \Xi^*[S_i],
\]
and thus
\[
H^*_\beta(\Xi^*[S_1,S_2,\dots,S_k]) = \bigotimes_i H^*_\beta(\Xi^*[S_i]).
\]
We therefore need only to determine $H^*_\beta(\Xi^*[\{1,2,\dots,n\}])$.

The main result of this section is the following.

\begin{theorem}\label{thm52}
There is a natural graded module isomorphism
\[
H^*_\beta(\Xi^*[\{1,2,\dots,2n+1\}])
\cong
\tilde\Lambda[\{1,2,\dots,2n+1\}]\otimes \FF_2,
\]
and thus the former is concentrated in degree $n$.  For even lengths,
\[
H^*_\beta(\Xi^*[\{1,2,\dots,2n\}])=0.
\]
\end{theorem}

\begin{proof}
We first observe that there is a map from $\tilde\Lambda(S)\otimes \FF_2$
to $H^*_\beta(\Xi^*(S))$, taking $\tau_{ijk}$ to $\Pi_{\{i,j,k\}}$.  
To see that this is well-defined,
we check by (\ref{sot}) that $\Pi_{\{i,j,k\}}$ is $\beta$-closed,
and verify that the relations
hold modulo the image of $\beta$.  Indeed, we find
\begin{align}
\Pi_{\{i,j,k\}}^2 &= 0 \notag\\
\Pi_{\{i,j,k\}}\Pi_{\{j,k,l\}} &= 
\Pi_{\{i,j,k\}}\Pi_{\{i,j,k,l\}}
+
\Pi_{\{j,k,l\}}\Pi_{\{i,j,k,l\}}
+
\Pi_{\{i,j,k,l\}}^2\notag\\
&=
\beta(\Pi_{\{i,j,k,l\}})\notag
\end{align}
by (\ref{hrt}) and (\ref{bw}).
For the 5-term relation, we take the identity
\[
\Pi_{\{i,j,k\}}\Pi_{\{k,l,m\}}
=
\Pi_{\{i,j,k\}}\Pi_{\{i,j,k,l,m\}}
+
\Pi_{\{k,l,m\}}\Pi_{\{i,j,k,l,m\}}
+
\beta(\Pi_{\{i,j,k,l,m\}})
\]
and sum over cyclic shifts to obtain
\begin{align}
\Pi_{\{i,j,k\}}\Pi_{\{k,l,m\}}&+
\Pi_{\{j,k,l\}}\Pi_{\{l,m,i\}}+
\Pi_{\{k,l,m\}}\Pi_{\{m,i,j\}}\notag\\
&+\Pi_{\{l,m,i\}}\Pi_{\{i,j,k\}}+
\Pi_{\{m,i,j\}}\Pi_{\{j,k,l\}}
=
 \beta(\Pi_{\{i,j,k,l,m\}}),\notag
\end{align}
since each term $\Pi_{\{i,j,k\}}\Pi_{\{i,j,k,l,m\}}$ appears exactly twice
in the sum.

We will prove the theorem by induction, using a spectral sequence.
Consider the filtration of $\Xi^*[\{1,\ldots,n\}]$ given by
\[
V_i=\Pi_{\{1,2,\dots,n\}}^i \Xi^*_n,
\]
for $1\le i\le n-2$; each of these is a homogeneous ideal of $\Xi^*_n$.  
The quotient $V_i/V_{i+1}$ is spanned by monomials $\prod \Pi_S^{d_S}$ having $d_{\{1,2,\ldots,n\}}=i$.
By Theorem \ref{thm45}, we therefore
have a natural identification (including the action of $\beta$)
\[
V_i/V_{i+1} = \Pi_{\{1,2,\ldots,n\}}^i\bigoplus_{\pi:|\pi|>i+1} \Xi^*[\pi]\,\cong\bigoplus_{\pi:|\pi|>i+1} \Xi^*[\pi].
\]
The filtration $V_i$ gives rise to a spectral sequence converging to
$H^*_\beta(\Xi^*_n)$ with first page
\[
E_1^{p,q} = H^{p+q}_\beta(V_p/V_{p+1})
          \cong \bigoplus_{\pi:|\pi|>p+1} H^q_\beta(\Xi^*[\pi]).
\]
Since each set in each partition $\pi$ has size $<n$, we find by the
inductive hypothesis that
\[
H^*_\beta(\Xi^*[\pi]) = 0
\]
unless all parts of $\pi$ have odd order, in which case it is the homology
of the interval $[\hat 0,\pi]$, where $\hat 0$ is the partition into
singletons.  
More precisely, the interval splits naturally as a product of
intervals $[\hat 0,\{1,2,\dots,k\}]$, and 
by Proposition \ref{Kun},
its homology is 
the corresponding tensor product.
This agrees
with the product decomposition coming from $\Xi^*$.  In particular, the
cohomology is concentrated in degree $(n-|\pi|)/2$.

Now, the filtration $V_p$ splits via the canonical monomials; from this
splitting, we conclude that the $d_1$ differential
$E_1^{p,q}\to E_1^{p+1,q}$ 
is induced by multiplication by $p\cdot \Pi_{\{1,2,\ldots,n\}}$.
This is zero for even $p$ and surjective for $p$ odd. 
We conclude that
$E_2^{2p,q} = 0$ and that
\[
E_2^{2p+1,q}
\cong
\bigoplus_{\pi:|\pi|=2p+3} H^q_\beta(\Xi^*[\pi]).
\]
In particular, if $n$ is even, $n-|\pi|$ is odd for all set partitions that
appear, and therefore $E_2^{p,q}=0$ for all $p,q$. It follows that $H_\beta^*(\Xi^*[\{1,2,\ldots,n\}])=0$ as desired.

If $n$ is odd, $E_2^{2p+1,q}=0$ unless $q=(n-2p-3)/2$, when by induction
\[
\begin{split}
E_2^{2p+1,(n-2p-3)/2}
&\cong
\bigoplus_{\pi:|\pi|=2p+3} \tilde\Lambda[\pi]\otimes\FF_2 \\
&\cong
\bigoplus_{\pi:|\pi|=2p+3} \tilde H_q([\hat{0},\pi])
=
W_q(\Piodd{n}).
\end{split}
\]
Here's an example of what the spectral sequence looks like for $n=9$:
\[
\put(-5,-73){\line(0,1){80}}
\put(-5,-73){\line(1,0){260}}
\xymatrix@R=.3cm@C=.4cm{
\scriptstyle E^{1,3}\ar@{}[r]|{=W_3}\ar[drr]&&&&&&\\
\scriptstyle E^{1,2}\ar[r]^{\scriptscriptstyle\sim}&\scriptstyle E^{2,2}&
\scriptstyle E^{3,2}\ar@{}[r]|{=W_2}\ar[drr]&&\ar@{}[urr]|{\textstyle 0}\\
\scriptstyle E^{1,1}\ar[r]^{\scriptscriptstyle\sim}&\scriptstyle E^{2,1}&
\scriptstyle E^{3,1}\ar[r]^{\scriptscriptstyle\sim}&\scriptstyle E^{4,1}&
\scriptstyle E^{5,1}\ar@{}[r]|{=W_1}\ar[drr]+<-20pt,4pt>&\\
\scriptstyle E^{1,0}\ar[r]^{\scriptscriptstyle\sim}&\scriptstyle E^{2,0}&
\scriptstyle E^{3,0}\ar[r]^{\scriptscriptstyle\sim}&\scriptstyle E^{4,0}&
\scriptstyle E^{5,0}\ar[r]^{\scriptscriptstyle\sim}&\scriptstyle E^{6,0}&
\scriptstyle E^{7,0}=W_0
}
\]

We claim that the $d_2$ differential $E_2^{2p+1,q}\to E_2^{2p+3,q-1}$ agrees with
the canonical differential on $\tilde\Lambda$ (recall Theorem \ref{thm34} and Lemma \ref{lem317}).
This would then imply that
\[
E_3^{*,*}=E_3^{1,(n-3)/2}=\tilde\Lambda[\{1,2,\ldots,n\}]\otimes\FF_2
\]
and the theorem would follow.

This differential again corresponds to multiplication by
$\Pi_{\{1,2,\dots,n\}}$: it takes an element $x\in\ker_\beta(V_{2p+1})$ with
$\Pi_{\{1,\dots,n\}}x\in V_{2p+3}$ to
\[
d_2\big(x+V_{2p+2}\big)=\Pi_{\{1,\dots,n\}}x+V_{2p+4}.
\]
We may as well assume that $x\in \Pi_{\{1,\dots,n\}}^{2p+1}\Xi^*[\pi]$
for some partition $\pi$ with $2p+3$ parts i.e. that 
\[
x=\big(\Pi_{\{1,\ldots,n\}}^{2p+1}\prod_i\Pi_{\pi_i}\big)y
\]
for some $y$. Multiplying by
$\Pi_{\{1,\dots,n\}}$, we use (\ref{dxd}) to compute
\[
d_2\big(x+V_{2p+2}\big)=
\big(\Pi_{\{1,\dots,n\}}^{2p+2}
\prod_i \Pi_{\pi_i}\big)y
\equiv
\big(\sum_i
\Pi_{\{1,\dots,n\}}^{2p+3}
\prod_{j\ne i} \Pi_{\pi_j}\big)y
\]
modulo $V_{2p+4}$.
In particular, the $d_2$ differential is induced by the natural
maps $\delta:=$``divide by $\Pi_{\pi_i}$'' on each piece of $\pi$
\[
\delta:H^*_\beta(V_1/V_2)\to H^{*-1}_\beta(V_0/V_1).
\]


So it suffices to prove that $\delta$ agrees with the canonical
differential in its action on $\tilde\Lambda$.  For this we proceed by
induction, observing that it holds for $\tilde\Lambda[\{1,2,3\}]$ (where
both maps take $\tau_{123}$ to 1), and that
it acts as a derivation on trees.
Indeed, using (\ref{hrt}), we check that
\[
\begin{split}
\delta\big((\Pi_Sy)(\Pi_Tz)\big)&=
\delta\big((\Pi_{S\cup T}\Pi_S+\Pi_{S\cup T}\Pi_T-\Pi_{S\cup T}^2)yz\big)\equiv\\
\delta\big(\Pi_{S\cup T}(\Pi_S+\Pi_T)yz\big)&=
(\Pi_S+\Pi_T)yz=
\delta(\Pi_Sy)\Pi_Tz+\delta(\Pi_Tz)\Pi_Sy.
\end{split}
\]
The theorem follows.
\end{proof}

Let $h^k(X,\Bbb Q)$ denote the $k$-th Betti number of a space $X$.

\begin{corollary} \label{th:bound}
We have the coefficient-wise upper bound
\[
\sum_{0\le n,k} h^k(M_{n+1}^\RR,\QQ) t^k \frac{u^n}{n!}
\le
\exp\big(\arcsin(u\sqrt{t})/\sqrt{t}\big).
\]
\end{corollary}

\begin{proof}
Indeed, we have
\[
h^k(M_{n+1}^\RR,\QQ)
\le
h^k_\beta(H^*(M_{n+1}^\RR,\FF_2))
=
h^k_\beta(\Xi_n^*)
=
\dim \tilde\Lambda(\{1,\ldots,n\}),
\]
and the latter has the stated exponential generating function by Corollary \ref{cor:hilser_t}.
\end{proof}

A twisted version may also be of interest.

\begin{corollary}\label{cor:twisted_Bockstein}
Let $\beta'$ be the twisted differential defined by
\[
\beta'(v) = \beta(v) + \Pi_{\{1,2,\dots,n\}} v.
\]
Then there is a canonical isomorphism
\[
H^*_{\beta'}(\Xi^*_{2n})
\cong
\tilde H_*(\Piodd{2n},\ZZ)\otimes \FF_2
\]
and
\[
H^*_{\beta'}(\Xi^*_{2n+1})=0.
\]
\end{corollary}


\section{The operad structure and the proof of the main results}
The main remaining step of the proof of the main theorem is to show that $
f_n^\ZZ $ is split-injective.  To do this we will use the operad structure
to construct elements of the homology $ H_*(M_n) $ which pair
upper-triangularly with the images under $ f_n $ of the basis vectors in $
\Lambda_n $.  Along the way we will also determine the structure of the
homology operad.

\subsection{The cyclic (co)operad $\Lambda$}

The collection of algebras $\Lambda_n$ forms 
a cyclic cooperad $\Lambda$. 
To describe this cyclic cooperad structure, it will be
convenient to replace the index set $\{1,2,\dots,n\}$ by an arbitrary
nonempty finite set.  Thus let $\Lambda\langle S\rangle$ be the
$\ZZ$-algebra with antisymmetric generators $\omega_{ijkl}$ for $i,j,k,l\in
S$ satisfying the relations \eqref{om1}, \eqref{om2}, \eqref{om3}.
Similarly, let $\Lambda(S)$ denote the $\ZZ$-algebra with antisymmetric
generators $\nu_{ijk}$ for $i,j,k\in S$ satisfying the relations of
$\Lambda_{|S|+1}$.

Given nonempty finite sets $S$, $T$ and a function $f:S\to T$, define a
function $h_t:S\to f^{-1}(t)\sqcup \{t\}$ for each $t\in T$ by
\[
h_t(s) = \begin{cases}
s, & \text{if }f(s)=t,\\
t, & \text{otherwise.}
\end{cases}
\]
Then we define a homomorphism (which will be the structure map of the
cooperad)
\[
\Delta_f:
\Lambda\langle S\rangle
\to
\Lambda\langle T\rangle\otimes
\bigotimes_{t\in T}
\Lambda\langle f^{-1}(t)\sqcup\{t\}\rangle
\]
in degree one as follows:
\[
\Delta_f(\omega_{ijkl})
=
\omega^0_{f(i)f(j)f(k)f(l)}
+
\sum_{t\in T} \omega^t_{h_t(i)h_t(j)h_t(k)h_t(l)},
\]
where each $\omega_{ijkl}=0$ if two indices agree.  (The superscripts are
used merely to distinguish the generators of $\Lambda\langle T\rangle$,
denoted $\omega^0_{ijkl}$, from the generators of each $\Lambda\langle
f^{-1}(t)\sqcup \{t\}\rangle$.)  Note that this implies that at most one
term appears on the right.  (N.b., since we are defining the map on
generators, the order of the tensor product is irrelevant.  Also, if
$|f^{-1}(t)|<3$, we may freely omit the corresponding factor on the right,
since $\Lambda\langle S\rangle=\ZZ$ if $1\le |S|\le 3$.)   This should be
viewed as corresponding to the geometric operad map
\[
M_T\times\prod_{t\in T} M_{f^{-1}(t)\sqcup\{t\}}
\to
M_S
\]
obtained by gluing together each pair of points labelled $t$.

The map $\Delta_f$ 
is easily seen to respect relation (\ref{om1}).
For relation (\ref{om2}),
we find that either one of the two generators maps to 0 or both map to
elements of the same algebra; in either case, the relation automatically
holds.  To check relation (\ref{om3}), there's a few more cases. Either all
the generators map to the same algebra, or at least one generator in each
monomial maps to 0, or we have
$f(i)\!=\!f(j)\!=\!f(k)\!\not=\!f(l)\!=\!f(m)\!=\!f(p)$, or we have
$f(i)\!=\!f(j)\!=\!f(k)$, $|\{f(i),f(l),f(m),f(p)\}|=4$ (up to cyclic
permutation of the indices).  In either of the latter cases, the relation
is easily verified.

Given another map $g:T\to U$, we define for each $u\in U$ a map
\[
f_u:(g\circ f)^{-1}(u)\sqcup \{u\}\to g^{-1}(u)\sqcup\{u\}
\]
by $f_u(s)=f(s)$ for $s\in (g\circ f)^{-1}(u)$, and $f_u(u)=u$.

\begin{theorem}
Given nonempty finite sets $S$, $T$, $U$ and functions $f:S\to T$, $g:T\to
U$, we have the identity
\[
(\Delta_g\otimes 1)\circ \Delta_f
=
(1\otimes \bigotimes_{u\in U}\Delta_{f_u})
\circ \Delta_{g\circ f}
\]
(up to the symmetry of the tensor product).

In other words, the maps $\Delta_f$ furnish the algebras $\Lambda\langle
S\rangle$ with the structure of a cyclic cooperad in the category of
superalgebras.
\end{theorem}

\begin{proof}
It suffices to check the relation on the generators of $\Lambda\langle S\rangle$.
And indeed, both
sides map $\omega_{ijkl}$ to
\[
\begin{split}
\omega^0_{g(f(i))g(f(j))g(f(k))g(f(l))}
&+
\sum_{u\in U} \omega^u_{h_u(f(i))h_u(f(j))h_u(f(k))h_u(f(l))}\\
&+
\sum_{t\in T} \omega^t_{h_t(i)h_t(j)h_t(k)h_t(l)}.
\end{split}
\]
\vspace{-.7cm}

\end{proof}

\noindent The following proposition follows easily from definitions. 

\begin{proposition}
The above cyclic cooperad structure is compatible with the natural
map $f_n^{\Bbb \ZZ_2}:
\Lambda_n\to H^*(M_n,\ZZ_2)/H^*(M_n,\ZZ_2)\langle 2\rangle $
and the cyclic cooperad structure induced on cohomology by the
cyclic operad structure of $M_n$.
\end{proposition}

We now turn to the $S_{n-1}$-symmetric algebras $\Lambda(S)$.  
If we add to $S$ a new label $\infty$, we know by Proposition \ref{nupres} that
$\nu_{ijk}\mapsto\omega_{ijk\infty}$ induces an isomorphism
$\Lambda(S)\cong\Lambda\langle S\sqcup\{\infty\}\rangle$.

Given a function $f:S\to T$, we can extend it to a function
$f:S\sqcup\{\infty\}\to T\sqcup\{\infty\}$ by taking $f(\infty)=\infty$.
This then gives a map which we abusively denote $\Delta_f$:
\[
\Delta_f:
\Lambda(S)
\to
\Lambda(T)
\otimes
\bigotimes_{t\in T}\Lambda(f^{-1}(t)).
\]
This map then satisfies the axioms of a (noncyclic) operad.  On the generators
$\nu$, we have
\[
\Delta_f(\nu_{abc})
=
\begin{cases}
\nu^0_{f(i)f(j)f(k)}, & \text{if }|f(\{i,j,k\})|=3,\\
0, & \text{if }|f(\{i,j,k\})|=2,\\
\nu_{ijk}^t, & \text{if }f(i)\!=\!f(j)\!=\!f(k)\!=:\!t.
\end{cases}
\]

Now, if we quotient by the augmentation ideal of $\Lambda(T)$, we obtain a
map $\eta_f: \Lambda(S)\to \bigotimes_{t\in T}\Lambda(f^{-1}(t))$.

\begin{proposition}
The map $\eta_f$ is a split surjection.
\end{proposition}

\begin{proof}
Indeed, for each $t\in T$, we have a map $\Lambda(f^{-1}(t))\to \Lambda(S)$
induced from the inclusion map; taking the product of these maps gives the
desired splitting.
\end{proof}

Since $\Lambda$ is a cooperad, its dual $\Lambda^*$ is an operad.
The simplest operation of $\Lambda^*$ not in the suboperad of commutative
superalgebras (i.e. the degree 0 part) is the ternary operation $\tau$, of
degree -1, corresponding to the linear functional on $\Lambda(\{1,2,3\})$
that takes 1 to 0 and $\nu_{123}$ to 1.  Together with the supercommutative
product $\cdot\,$, the operation $\tau$ satisfies the following relations:
\begin{align}
&\tau(x,y,z) = (-1)^{1+|x||y|} \tau(y,x,z)= (-1)^{1+|y||z|} \tau(x,z,y)\label{taurel1}\\
&\tau(w,x,y\cdot z) = \tau(w,x,y)\cdot z+(-1)^{|y||z|}\tau(w,x,z)\cdot y,\label{taurel2}
\end{align}
and a 10-term relation stating that the various permutations
(with appropriate signs compatible with superantisymmetry) of
$\tau(\tau(v,w,x),y,z)$ sum to 0.  (That some 10-term relation holds
follows from the fact that $\Lambda_{5+1}[2]$ is free of rank $9$, while
the space of compositions of $\tau$ is free of rank $10$; that it has the
stated form follows from the fact that $\Hom(\Lambda_{5+1}[2],\ZZ)$ does
not contain a copy of the sign representation of $S_5$.)  

Let $\Lambdasharp$ be the operad generated by $\cdot$ and $\tau$, with
relations given by (\ref{taurel1}), (\ref{taurel2}) and the above 10-term
relation.  We have a map $ \Lambdasharp \rightarrow \Lambda^* $ which we
will soon show is an isomorphism.

If we ignore the product, the operad $ {\Lambdasharp}'(\bullet) $ generated
by $\tau$ with superantisymmetry and the 10-term relation is a twisted
version of the Lie $2$-algebra operad $ HW(\bullet) $ discussed in
\cite{HW}.  One of their main results is that $HW(n)$ is isomorphic to
$H^*(\Piodd{n})$.  It follows from Theorem \ref{th:dimposethom} that the
exponential generating function for $\dim(HW(n))$ is
$\arcsin(u\sqrt{t})/\sqrt{t}$.  Since the twisting just amounts to
tensoring with the sign operad, the same generating function applies to our
twisted operad.

Using (\ref{taurel2}), any operation in $\Lambdasharp$ can be expressed as
a sum of operations where $\cdot$ is only used after $\tau$.  It follows
that the exponential generating function for the number of additive
generators of $\Lambdasharp$ is bounded above by
$\exp(\arcsin(u\sqrt{t})/\sqrt{t})$.  Recall that by Corollary
\ref{cor:hilser_t} this was also the exponential generating function for
the ranks of $ \tilde\Lambda_n $.  A similar result holds for $ \Lambda_n$.

\begin{thm} \label{th:basis}
The basic forests form a basis of $\Lambda(S)$, which is thus a free
$\ZZ$-module.  The corresponding exponential generating function is
\[
P(u,t)=\exp\big(\arcsin(u\sqrt{t})/\sqrt{t}\big),
\]
and the Hilbert series of $\Lambda_n$ is
\[
P_n(t)=\prod_{0\le k<(n-3)/2}(1+(n-3-2k)^2t).
\]
\end{thm}

\begin{proof}
The argument in Proposition \ref{prop312} was purely combinatorial, and
thus changing signs in the relations will have no effect.  As a result,
$\Lambda(S)$ is also spanned by basic forests.

Unfortunately, by twisting the $S_n$ action, we have destroyed the
canonical isomorphism with partition homology, making the argument in
Theorem \ref{thm:homology_algebra} no longer valid.  Recall that the main
tool of Theorem \ref{thm:homology_algebra} was a pairing between saturated
chains of odd partitions and triangle trees.  We can achieve the same
effect using ternary forests and the operad structure.

To each saturated chain of odd partitions, we associate a ternary forest,
with a node for each set that appears as the part of some partition in the
chain.  To each ternary tree $ U $ , let $ \tau_U \in \Lambda^*(S) $ denote
the result of composition of the operation $\tau $ according to the tree $
U $; we then extend this to forests using the product operation $\cdot$.
Recall that to each triangle forest $ F $, we have a monomial $ \nu^F \in
\Lambda(S) $.  Now given $ F $ a triangle forest and $ G $ a ternary forest, we
can consider the pairing $ \langle \tau_G, \nu^F \rangle $.

The action of composition with $\tau$ (corresponding to a partition into
three components) is such that a triangle forest has nonzero image iff
there is a triangle in the forest that hits each component, while all other
triangles are contained in a component; the action of composition with
$\cdot$ is such that a triangle forest has nonzero image iff each component
of the forest is contained in a component of the composition.  In
particular, the pairing between ternary forests and triangle forests takes on
only the values $\pm 1$ and $0$, nonzero iff there exists an ordering on
the triangles inducing a partition chain with the given associated ternary
forest.

\begin{minipage}{12.5cm}
\vspace{.3cm}
\centerline{
\psfig{file=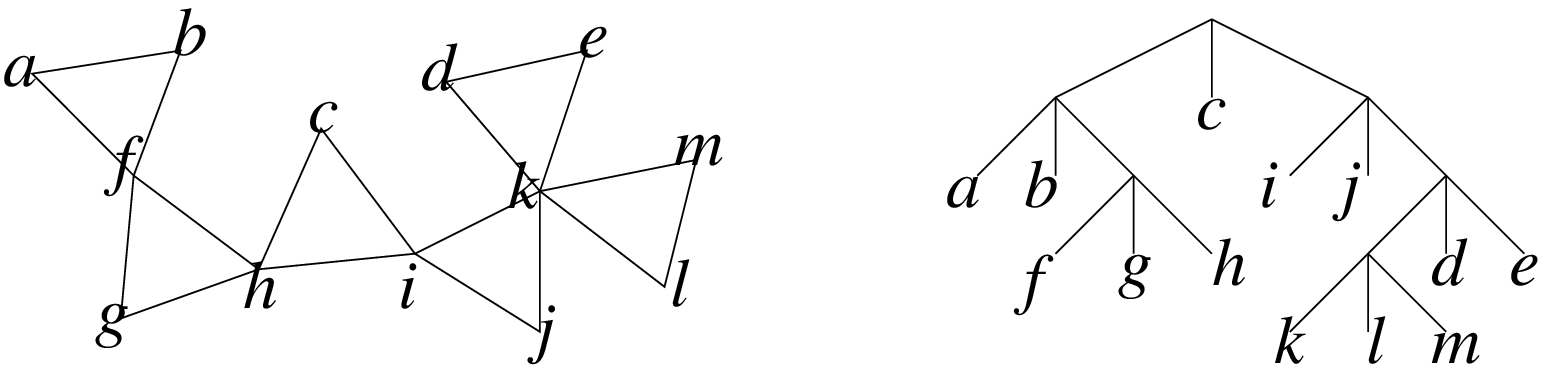,width=8cm}
}\vspace{.1cm}

\begin{center}
{\it
A triangle tree and a ternary tree that pair non-trivially.
}
\end{center}
\vspace{.3cm}
\end{minipage}

By mimicking the proof of Theorem \ref{thm:homology_algebra}, with ternary
forests instead of chains of odd partitions we obtain the result; note that
the action of $\cdot$ is such that the pairing between forests with the
same component partition is just (up to sign) the product of the pairings
between the individual trees.
\end{proof}

The proof of the theorem shows that the pairing between $\Lambdasharp$ and
$\Lambda$ induces a surjective map $\Lambdasharp\to\Hom(\Lambda,\ZZ)$.
Moreover, the number of additive generators of $\Lambdasharp$ is bounded
above by the rank of $\Lambda$.  It follows that $\Lambdasharp$ is torsion
free and that the natural map $\Lambdasharp\to\Lambda^*$ is an isomorphism.
We conclude the following.

\begin{corollary}
The operad $\Lambda^*$ is presented over the operad of commutative
superalgebras by a ternary operation $\tau$ satisfying 
(\ref{taurel1}), (\ref{taurel2}) and the 10-term relation $\text{\rm Alt}(\tau\circ(\tau\otimes\text{\rm Id}\otimes\text{\rm Id}))=0$.
\end{corollary}

\subsection{Determination of the cohomology ring}
We are now in a position to prove our main theorems.
\begin{thm}\label{cor:f_injective}
The map $f_n:\Lambda_n\to H^*(M_n,\ZZ)/H^*(M_n,\ZZ)\langle 2\rangle $ is
injective and splits as a map of $\ZZ$-modules.
\end{thm}

\begin{proof}
Recall that $ \Lambda_n $ is a free $\ZZ $ module with basis of monomials $
\nu^F $ indexed by basic triangle forests $ F $.  It suffices to find a
collection of elements in $ H_*(M_n, \ZZ) $, also indexed by basic triangle
forests which pair upper triangularily with the images $ f_n(\nu^F) $.

To do this we follow the proof of Theorem \ref{thm:homology_algebra}, or
more precisely its modification in Theorem \ref{th:basis}.  In the proof of
Theorem \ref{thm:homology_algebra}, for each basic triangle forest (more
precisely, tree, but by the product structure, this extends to forests), we
associated a saturated chain of set partitions, which paired with the
desired property.  In the proof of Theorem \ref{th:basis}, we used this
saturated chain of set partitions to build a ternary forest.

For each ternary tree $ U $, we can consider a map $g_U : M_4^d\to M_n$
which is given by gluing curves according to the tree $ U $.  Let $ \rho^U
$ denote the image under the map $ g_U $ of the fundamental class in $
H_*(M_n) $.  These elements $ \rho^U $ are exactly what would be obtained
by using $ U $ to specify a particular composition of the homology operad
generator $ \tau \in H_1(M_4, \QQ) $ (after tensoring with $ \QQ $).  More
generally, for a ternary forest $G$, we consider the associated map
$g_G:M_4^{d'} \times M_c\to M_n$, where $c$ is the number of components of
$G$, first gluing together curves according to the component trees of the
forest, then gluing their roots to $M_c$.  The class $\rho^G$ is then the
image of the product of the fundamental class of $M_4^d$ with the point
class in $M_c$, and again corresponds to the appropriate element of the

homology operad.

The pairing between $ \rho^G $ and $ f_n(\nu^F) $ is again described
combinatorially as in the proof of Theorem \ref{th:basis}, and thus we
again have the desired upper triangularity property.

In particular, we have constructed a collection of linear functionals on
$H^*(M_n,\ZZ)/H^*(M_n,\ZZ)\langle 2\rangle $ that, when evaluated on our
basis of $\Lambda_n$, produce a triangular matrix with unit diagonal.
Composing these functionals with the inverse of that upper triangular
matrix produces a left inverse of $f_n$.  It follows that $\Lambda_n$ is a
direct summand of $H^*(M_n,\ZZ)/H^*(M_n,\ZZ)\langle 2\rangle $.
\end{proof}

So we can write $H^*(M_n,\ZZ)/H^*(M_n,\ZZ)\langle 2\rangle = \Lambda_n
\oplus B $.  At this point we do not know much about the complement $ B $.
Here we need the information from section \ref{se:bound}.

\begin{thm}
$ B $ is a product of odd torsion modules and $ H^*(M_n, \ZZ) $ has no
  4-torsion.
\end{thm}

\begin{proof}
Let us split $ H^*(M_n, \ZZ) $ (non-canonically) as $ H^*(M_n, \ZZ) = F
\oplus G \oplus H \oplus I $, where $ F $ is free, $ G $ is a product of
odd torsion modules, $ H $ is a product of $ \ZZ/2\ZZ $ and $ I $ is a
product of $ \ZZ/2^k\ZZ $ for $k \ge 2 $.  Then, $
H^*(M_n,\ZZ)/H^*(M_n,\ZZ)\langle 2\rangle $ is isomorphic to $ F \oplus G
\oplus I' $ where $ I' = I / I\langle 2\rangle $.  By Theorem
\ref{cor:f_injective}, we see that the $ \mathrm{rk}(F) \ge
\mathrm{rk}(\Lambda_n) $.

On the other hand, $H^*(M_n,\ZZ/4\ZZ)/H^*(M_n,\ZZ/4\ZZ)\langle 2\rangle $ is an $\FF_2 $ vector space of dimension $ \mathrm{rk}(F) $ plus twice the number of factors in $ I $.

We know from Theorem \ref{thm52} that there are isomorphisms
\begin{equation*}
\tilde\Lambda_n\otimes \FF_2 \cong \Lambda_n\otimes \FF_2 \to
H^*(M_n,\ZZ/4\ZZ)/H^*(M_n,\ZZ/4\ZZ)\langle 2\rangle.
\end{equation*}

Since $ \Lambda_n $ is free, the dimension of $ \Lambda_n \otimes \FF_2 $ is $ \mathrm{rk}(\Lambda_n) $.  Combining together these observations we see that
\begin{equation*}
\begin{split}
\mathrm{rk}(\Lambda_n) &= \mathrm{dim}(H^*(M_n,\ZZ/4\ZZ)/H^*(M_n,\ZZ/4\ZZ)\langle 2\rangle) \\
&= \mathrm{rk}(F) +  2(\text{\# of factors in } I) \ge \mathrm{rk}(F) \ge \mathrm{rk}(\Lambda_n).
\end{split}
\end{equation*}
Hence we conclude that $ \mathrm{rk}(F) = \mathrm{rk}(\Lambda_n) $ and that $ I = 0 $.  Since $ I = 0 $, there is no 4-torsion.  Also we see that the complement $ B $ is a product of odd torsion modules.
\end{proof}

\begin{corollary}
The map $ f^{\ZZ_2} : \Lambda_n \otimes \ZZ_2 \rightarrow H^*(M_n, \ZZ_2)/H^*(M_n, \ZZ_2)\langle 2 \rangle $ is an isomorphism.
\end{corollary}

As a corollary of our proof we also see the following.
\begin{corollary}
$ H_*(M_n, \QQ) $ has a basis given by fundamental classes tensor point classes coming
  from maps $ M_4^d \times M_c\rightarrow M_n $ corresponding to basic
  triangle forests.
\end{corollary}

\end{document}